\documentclass[11pt,reqno]{amsart}
\usepackage{amsthm}
\usepackage{amssymb}
\usepackage{latexsym}
\usepackage{multicol}
\usepackage{verbatim,enumerate}
\usepackage{amscd}

\usepackage[usenames]{color}
\usepackage{hyperref}
\usepackage{amsmath, amscd}

\advance\textwidth by 1.2in \advance\oddsidemargin by -.6in \advance\evensidemargin by -.6in
\newtheorem*{cor}{Corollary}%[section]
\newtheorem*{lem}{Lemma}
\newtheorem*{prop}{Proposition}

\theoremstyle{definition}

\theoremstyle{definition}
\newtheorem{thm}{Theorem}
\newtheorem*{thm*}{Theorem}

\newtheorem*{rem}{Remark}

\newenvironment{pf}{\proof}{\endproof}
\newcounter{cnt}
\newenvironment{enumerit}{\begin{list}{{\hfill\rm(\roman{cnt})\hfill}}{%
\settowidth{\labelwidth}{{\rm(iv)}}\leftmargin=\labelwidth%
\advance\leftmargin by \labelsep\rightmargin=0pt\usecounter{cnt}}}{\end{list}} \makeatletter
\def\mydggeometry{\makeatletter\dg@YGRID=1\dg@XGRID=20\unitlength=0.003pt\makeatother}
\makeatother \theoremstyle{remark}

\numberwithin{equation}{section}

 \DeclareMathOperator{\Ht}{ht} 
\let\bwdg\bigwedge
\def\bigwedge{{\textstyle\bwdg}}

\begin{document}

\newcommand{\thmref}[1]{Theorem~\ref{#1}}
\newcommand{\secref}[1]{Section~\ref{#1}}
\newcommand{\lemref}[1]{Lemma~\ref{#1}}
\newcommand{\propref}[1]{Proposition~\ref{#1}}
\newcommand{\corref}[1]{Corollary~\ref{#1}}
\newcommand{\remref}[1]{Remark~\ref{#1}}
\newcommand{\defref}[1]{Definition~\ref{#1}}
\newcommand{\er}[1]{(\ref{#1})}
\newcommand{\id}{\operatorname{id}}
\newcommand{\ord}{\operatorname{\emph{ord}}}
\newcommand{\sgn}{\operatorname{sgn}}
\newcommand{\wt}{\operatorname{wt}}
\newcommand{\tensor}{\otimes}
\newcommand{\from}{\leftarrow}
\newcommand{\nc}{\newcommand}
\newcommand{\rnc}{\renewcommand}
\newcommand{\dist}{\operatorname{dist}}
\newcommand{\qbinom}[2]{\genfrac[]{0pt}0{#1}{#2}}
\nc{\cal}{\mathcal} \nc{\goth}{\mathfrak} \rnc{\bold}{\mathbf}
\renewcommand{\frak}{\mathfrak}
\newcommand{\supp}{\operatorname{supp}}
\newcommand{\Irr}{\operatorname{Irr}}
\newcommand{\psym}{\mathcal{P}^+_{K,n}}
\newcommand{\psyml}{\mathcal{P}^+_{K,\lambda}}
\newcommand{\psymt}{\mathcal{P}^+_{2,\lambda}}
\renewcommand{\Bbb}{\mathbb}
\nc\bomega{{\mbox{\boldmath $\omega$}}} \nc\bpsi{{\mbox{\boldmath $\Psi$}}}
 \nc\balpha{{\mbox{\boldmath $\alpha$}}}
 \nc\bpi{{\mbox{\boldmath $\pi$}}}
  \nc\bxi{{\mbox{\boldmath $\xi$}}}
\nc\bmu{{\mbox{\boldmath $\mu$}}} \nc\bcN{{\mbox{\boldmath $\cal{N}$}}} \nc\bcm{{\mbox{\boldmath $\cal{M}$}}} \nc\blambda{{\mbox{\boldmath
$\lambda$}}}\nc\bnu{{\mbox{\boldmath $\nu$}}}

\newcommand{\Tmn}{\bold{T}_{\lambda^1, \lambda^2}^{\nu}}

\newcommand{\lie}[1]{\mathfrak{#1}}
\newcommand{\ol}[1]{\overline{#1}}
\makeatletter
\def\section{\def\@secnumfont{\mdseries}\@startsection{section}{1}%
  \z@{.7\linespacing\@plus\linespacing}{.5\linespacing}%
  {\normalfont\scshape\centering}}
\def\subsection{\def\@secnumfont{\bfseries}\@startsection{subsection}{2}%
  {\parindent}{.5\linespacing\@plus.7\linespacing}{-.5em}%
  {\normalfont\bfseries}}
\makeatother
\def\subl#1{\subsection{}\label{#1}}
 \nc{\Hom}{\operatorname{Hom}}
  \nc{\mode}{\operatorname{mod}}
\nc{\End}{\operatorname{End}} \nc{\wh}[1]{\widehat{#1}} \nc{\Ext}{\operatorname{Ext}} \nc{\ch}{\text{ch}} \nc{\ev}{\operatorname{ev}}
\nc{\Ob}{\operatorname{Ob}} \nc{\soc}{\operatorname{soc}} \nc{\rad}{\operatorname{rad}} \nc{\head}{\operatorname{head}}
\def\Im{\operatorname{Im}}
\def\gr{\operatorname{gr}}
\def\mult{\operatorname{mult}}
\def\Max{\operatorname{Max}}
\def\ann{\operatorname{Ann}}
\def\sym{\operatorname{sym}}
\def\loc{\operatorname{loc}}
\def\Res{\operatorname{\br^\lambda_A}}
\def\und{\underline}
\def\Lietg{$A_k(\lie{g})(\bsigma,r)$}

 \nc{\Cal}{\cal} \nc{\Xp}[1]{X^+(#1)} \nc{\Xm}[1]{X^-(#1)}
\nc{\on}{\operatorname} \nc{\Z}{{\bold Z}} \nc{\J}{{\cal J}} \nc{\C}{{\bold C}} \nc{\Q}{{\bold Q}}
\renewcommand{\P}{{\cal P}}
\nc{\N}{{\Bbb N}} \nc\boa{\bold a} \nc\bob{\bold b} \nc\boc{\bold c} \nc\bod{\bold d} \nc\boe{\bold e} \nc\bof{\bold f} \nc\bog{\bold g}
\nc\boh{\bold h} \nc\boi{\bold i} \nc\boj{\bold j} \nc\bok{\bold k} \nc\bol{\bold l} \nc\bom{\bold m} \nc\bon{\bold n} \nc\boo{\bold o}
\nc\bop{\bold p} \nc\boq{\bold q} \nc\bor{\bold r} \nc\bos{\bold s} \nc\boT{\bold t} \nc\boF{\bold F} \nc\bou{\bold u} \nc\bov{\bold v}
\nc\bow{\bold w} \nc\boz{\bold z} \nc\boy{\bold y} \nc\ba{\bold A} \nc\bb{\bold B} \nc\bc{\bold C} \nc\bd{\bold D} \nc\be{\bold E} \nc\bg{\bold
G} \nc\bh{\bold H} \nc\bi{\bold I} \nc\bj{\bold J} \nc\bk{\bold K} \nc\bl{\bold L} \nc\bm{\bold M} \nc\bn{\bold N} \nc\bo{\bold O} \nc\bp{\bold
P} \nc\bq{\bold Q} \nc\br{\bold R} \nc\bs{\bold S} \nc\bt{\bold T} \nc\bu{\bold U} \nc\bv{\bold V} \nc\bw{\bold W} \nc\bz{\bold Z} \nc\bx{\bold
x} \nc\KR{\bold{KR}} \nc\rk{\bold{rk}} \nc\het{\text{ht }}

\nc\toa{\tilde a} \nc\tob{\tilde b} \nc\toc{\tilde c} \nc\tod{\tilde d} \nc\toe{\tilde e} \nc\tof{\tilde f} \nc\tog{\tilde g} \nc\toh{\tilde h}
\nc\toi{\tilde i} \nc\toj{\tilde j} \nc\tok{\tilde k} \nc\tol{\tilde l} \nc\tom{\tilde m} \nc\ton{\tilde n} \nc\too{\tilde o} \nc\toq{\tilde q}
\nc\tor{\tilde r} \nc\tos{\tilde s} \nc\toT{\tilde t} \nc\tou{\tilde u} \nc\tov{\tilde v} \nc\tow{\tilde w} \nc\toz{\tilde z} \nc\woi{w_{\omega_i}}
\nc\chara{\operatorname{Char}}
\title{Demazure modules, Fusion products and $Q$--systems}
\author[Chari and Venkatesh]{Vyjayanthi Chari and R.Venkatesh}\address{\noindent Department of Mathematics, University of California, Riverside, CA 92521}
\email{chari@math.ucr.edu}
\address{\noindent Department of Mathematics, Tata Institute of Fundamental Research, Homi Bhabha Road, Mumbai 400005, India.}
\email{rvenkat@math.tifr.res.in}
\thanks{V.C. was partially supported by DMS-0901253 and DMS- 1303052.}
\begin{abstract} In this paper, we introduce a family of indecomposable finite--dimensional graded modules for the current algebra associated to a simple Lie algebra. These modules are indexed by an  $|R^+|$--tuple of  partitions $\bxi=(\xi^\alpha)$, where $\alpha$ varies over a  set $R^+$ of positive roots of $\lie g$ and we assume that they satisfy a natural compatibility condition. In the case when the $\xi^\alpha$ are all rectangular, for instance, we prove that these modules are  Demazure modules in various levels. As  a consequence we see that the defining relations of Demazure modules can be greatly simplified. We  use this simplified presentation to relate our results to the fusion products,  defined in \cite{FL}, of representations of the current algebra. We prove that the $Q$--system of \cite{HKOTY} extends to a canonical short exact sequence of  fusion products of representations associated to certain special partitions $\bxi$. Finally, in the last section we deal with the case of $\lie{sl}_2$ and prove that the modules we define are just fusion products of irreducible representations of the associated current algebra and give monomial bases for these modules.
\end{abstract}
\maketitle

\section*{Introduction} The current algebra associated to a simple Lie algebra $\lie g$  is just the Lie algebra of polynomial maps  $\bc\to \lie g$.  As a vector space it is isomorphic to  $\lie g\otimes \bc[t]$, where $\bc[t]$ is the polynomial ring in the indeterminate $t$ and the Lie bracket is given in the obvious way. Both $\lie g[t]$ and its universal enveloping algebra inherit a  grading coming from the natural grading on $\bc[t]$.  The study of the category of  graded  finite--dimensional representations of the  current algebra has been of interest in recent years for a variety of reasons. The work of \cite{KodNaoi}  relates  graded characters of certain representations to the Poincare polynomials of quiver varieties. The homological properties of the category are similar to those of the BGG category $\cal O$ for the simple Lie algebra. This similarity has been explored to some extent  in \cite{BBCKL} and \cite{BCM} and leads to connections with symmetric functions and Macdonald polynomials.

The subject also has many  connections with  problems arising in mathematical physics, for instance  the $X=M$ conjectures, see \cite{eddy}, \cite{deFK}, \cite{Naoi1}. To explain this further, recall that
for any simple Lie algebra  there exists an
associated system of equations known as a $Q$-system. These were first introduced by Kirillov and Reshetikhin \cite{KR1} as a  combinatorial tool to address the
question of completeness of the Bethe ansatz states in the diagonalization of the Heisenberg spin
chain based on an arbitrary Lie algebra. The Q-systems can be considered as a type of discrete dynamical system, and the
 equations   are indexed by  pairs $(\alpha, m)$ where $\alpha$ varies over the set of simple roots of $\lie g$ with respect to a fixed Cartan subalgebra $\lie h$ and $m$ varies over the non--negative integers. The equations are defined recursively in $m$ and have the remarkable property that their solutions are polynomials, given
appropriate initial data. The Q-system with special (singular) initial conditions was originally introduced \cite{KR1} as the
recursion relation satisfied by the characters of special finite-dimensional modules of the Yangian $Y(\lie g)$, the so-called Kirillov--Reshetikhin modules. In the case $\lie g = A_n,$ the same
recursion relation also occurs in other contexts: for instance in the study of   Toda flows in Poisson geometry \cite{GSV} and in
preprojective algebras \cite{GLSP}. In \cite{deFK}  the
combinatorial Kirillov-Reshetikhin conjecture  which are the completeness conjectures for
the generalized Heisenberg spin chains is established.

It was shown in  \cite{H1}, \cite{NakT1}, that the $Q$--systems are more than a formal equality of characters of representations of quantum affine algebras. To be precise, it is shown that the solutions to the $Q$--systems come from a family of short exact sequences of tensor products of suitable Kirillov--Reshetikhin modules. These give rise to a more precise equality of characters and the corresponding system of equations are called the $T$-systems. These systems have played an important role in the theory of cluster algebras through the work of \cite{HL}. In recent work \cite{MY}, the notion of $T$--systems was extended beyond Kirillov--Reshetikhin modules to include a much bigger family of modules and it is expected that many of the applications of $T$--systems will go over to the extended T--systems.

A different approach was adopted in \cite{Ckir} and developed further in \cite{CMkir}. The goal of these papers was to understand the  $q=1$ limit of the  solutions of $Q$--systems. It was shown that the solutions  are characters of certain reducible but indecomposable representations of the classical affine Lie algebra and its standard maximal parabolic subalgebra, the current algebra. Adopting this point of view has led to many further combinatorial connections with the theory of Macdonald polynomials \cite{BBCKL}, \cite{CI}. The more general philosophy behind this approach is that
 many interesting families of irreducible representations  of the  quantum affine algebra (associated to $\lie g$), when specialized to $q=1$ give  indecomposable representations of the current algebra (see \cite{CPweyl} for instance). These families  include the  Kirillov--Reshetikhin modules, the minimal affinizations  and the more general modules studied in \cite{MY}.
 The defining relations  of the Kirillov--Reshetikhin modules for  the current algebra was given in \cite{CMkir} and it is clear from the definition that they are graded by the non--negative integers. The graded $\lie g$--module decomposition was also given in this paper, and it is remarkable that the Hilbert series of this decomposition coincides with the deformed character formulae for the Kirillov--Reshetikhin modules  given in \cite{HKOTY} where the powers of $q$ appear for entirely different reasons from the study of solvable lattice models. As a consequence of this approach, it was proved in
 \cite{CMkir} and also in \cite{FoL}, that this version of the  Kirillov--Reshetikhin modules are isomorphic to certain Demazure modules in positive level representations of the affine Lie algebra.

Our paper is motivated by two questions. It is not hard to see that if $V_1$ and $V_2$ are representations of the quantum affine algebra, then the specialization of $V_1\otimes V_2$ to $q=1$ is not isomorphic to the tensor product of the specializations.  Indeed, this problem was  studied in \cite{CPweyl} and led to the definition of the local Weyl modules  for the  affine Lie algebra  and hence also the current algebra associated to $\lie g$. At about the same time, Feigin and Loktev introduced the notion of  a fusion product of graded representations of the current algebra. It was later proved in \cite{CL}, \cite{FoL}, \cite{Naoi} that a local Weyl module is a fusion product of  fundamental local Weyl modules.  Together with results established in \cite{BN}, \cite{Cbraid}, it follows that the local Weyl module is the $q=1$ limit of a tensor product of irreducible representations of quantum affine algebra.
In other words, at least in this very special case, it is true, that the specialization of  a tensor product of representation of the  quantum affine algebras is isomorphic to  the fusion product of representations of the current algebra. Sections 4  and 5 show that this remains true for certain kinds of Kirillov--Reshetikhin modules. In fact, the results of Section 5  can be interpreted as follows: they show that the short exact sequences (T--systems) studied in \cite{H1}, \cite{NakT1}, \cite{NakT2} of tensor products of Kirillov--Reshethikhin modules, give rise by specializing, to short exact sequences of fusion products of the Kirillov--Reshethikhin modules of the current algebra. We also anticipate that our approach should work in the more general situation studied in \cite{MY} and should provide insight into the $\lie g$--module structure of the classical limits of the corresponding large classes of modules considered in \cite{MY}.

The second motivation for our paper was to understand the presentation given in  \cite{FF} of  the fusion product of irreducible finite--dimensional representations of the current algebra associated to $\lie{sl}_2$. Although the definition of fusion products are given in a purely representation theoretic way, the proofs in \cite{FF} do not use this definition in any direct way.  Moreover, their relations do not  include (at least in any straightforward way), certain obvious ones which hold in the fusion product. In fact it is exactly these obvious relations which are critical for  our study of Demazure modules and $Q$--systems for arbitrary simple Lie algebras. In Section 6 of our paper, we recover the presentation of \cite{FF}  and more. We give three presentations of the fusion product and in one of them (see Proposition \ref{third}) the Demazure type relations show up in a natural way. Finally, we give a monomial bases for fusion products which are consistent with certain canonical short exact sequences of modules for $\lie{sl}_2[t]$.

The paper is organized as follows. Section one establishes the basic notation and elementary results needed in the rest of the paper. In Section two, we define the modules $V(\bxi)$, where $\bxi$ is a tuple of partitions indexed by the positive roots.  We give three equivalent presentations of these modules. In Section 3, we consider the case when   $\bxi$  consists of partitions which are either rectangular or certain types of fat hooks. We show that in this case, the module $V(\bxi)$ is  isomorphic to a  Demazure module of level $\ell$ and that the defining relations of these modules can be  highly simplified. These  are generalizations of results  which were   previously known only for the level one Demazure modules (see \cite{FoL}, \cite{Naoi}). Even in the level one case,  our methods are very different and our approach is more uniform. In Section 4, we recall the definition of fusion products and prove a \lq\lq Schur positivity" type result for the fusion product of Demazure modules (see Proposition 4.5). In Section 5, we show that there exists a short exact sequence of graded $\lie g[t]$--modules corresponding to the $Q$--system defined in \cite{HKOTY}. Finally, in Section 6, we study the case of $\lie{sl}_2$ and compute the dimension and give a monomial basis of the modules.

{\em Acknowledgements: Part of this work was done when the first author was visiting ICERM for the semester on Automorphic Forms, Combinatorial Representation Theory and Multiple Dirichlet Series. She thanks the organizers of the semester for this opportunity and also ICERM for the superb working conditions. The second author thanks Professor V.S. Sunder, the National Board of Higher Mathematics, India and the Department of Mathematics at UCR for the financial support which made this collaboration possible. }

\section{Preliminaries } Throughout this paper, $\bc$ will denote the field of complex numbers and $\bz$ (resp. $\bz_+$, $\bn$) the set of integers (resp. non--negative, positive integers). We let $\bc[t]$  be the polynomial ring in
an indeterminate $t$.

\subsection{}  Given a complex Lie algebra $\lie a$, we let $\bu(\lie a)$ be the corresponding universal enveloping algebra. The associated current Lie algebra is denoted by $\lie a[t]$: as a vector space it is just $\lie a\otimes \bc[t]$ and the Lie bracket is given in the natural way: $[a\otimes f, b\otimes g]=[a,b]\otimes f g$, for all $a,b\in\lie a$ and $f,g\in\bc[t]$.  Let $\lie a[t]_+$ be the ideal $\lie a\otimes t\bc[t]$. We shall freely identify $\lie a$ with the Lie subalgebra $\lie a\otimes 1$ of $\lie a[t]$ and we clearly have an isomorphism of vector spaces $$\lie a[t]=\lie a[t]_+\bigoplus\lie a, \qquad \bu(\lie a[t])\cong \bu(\lie a[t]_+)\otimes\bu(\lie a).$$  The degree grading on $\bc[t]$ defines a natural $\bz_+$--grading on $\lie a[t]$ and hence also on  $\bu(\lie a[t])$: an element of the form $(a_1\otimes t^{r_1})\cdots (a_s\otimes t^{r_s})$ has  grade $r_1+\cdots+r_s$. Denote by $\bu(\lie a[t])[r]$ the subspace of grade $r$.

\subsection{} A {\em graded} representation of $\lie a[t]$ is a $\bz$--graded vector space which admits a compatible Lie algebra action of $\lie a[t]$, i.e.,
$$V=\bigoplus_{r\in\bz}V[r],\qquad(\lie a\otimes t^s)V[r]\subset V[r+s],\ \  r\in\bz, \ \ s\in\bz_+.$$  If $V$ and $V'$ are graded $\lie a[t]$--modules, we say that $\pi:V\to V'$ is a morphism of graded $\lie a[t]$--modules if $\pi$ is a degree zero morphism of $\lie a[t]$--modules. For  $r\in\bz$,  let $\tau_r$  be the grading shift operator: if $V$ is a graded $\lie a[t]$--module then $\tau_rV$ is the graded $\lie a[t]$--module with the graded pieces shifted uniformly by $r$ and   the action of $\lie{a}[t]$ unchanged. If $M$ is an $\lie a$--module
and $z\in\bc$,  define a $\lie a[t]$--module structure on $M$ by: $(a\otimes t^r)m=z^ram$. We denote this module as $\ev_z M$. Clearly $\ev_zM$ is an irreducible $\lie a[t]$--module iff $M$ is an irreducible $\lie a$--module. Moreover the module $\ev_0M$ is a graded $\lie a[t]$--module and $\left(\ev_0M\right)[0]=M$. In particular, $\lie a[t]_+(\ev_0M) =0$.

\subsection{} From now on,  $\lie g$ will be an arbitrary simple finite--dimensional complex Lie algebra and  $\lie h$ will be   a fixed Cartan subalgebra of $\lie g$ and we assume that  $\dim\lie h=n$ or equivalently that the rank of $\lie g$ is $n$. Let $R$ be the set of roots of $\lie g$ with respect to $\lie h$.
The restriction of the Killing form of $\lie g$ to $\lie h$ induces an isomorphism between $\lie h$ and $\lie h^*$ and hence also a  symmetric non--degenerate form  $(\ , \ )$  on $\lie h^*$. We shall assume that this form on $\lie h^*$   is normalized so that the square length of a  long root is two. For $\alpha\in R$, let $t_\alpha\in\lie h$ be the element that maps to $\alpha$ and set, $$ d_\alpha= \frac{2}{(\alpha,\alpha)},\ \ \qquad h_\alpha=d_\alpha t_\alpha.$$

\vskip 6pt

 \noindent Let   $I=\{1,\cdots ,n\}$  and fix a set $\{\alpha_i: i\in I\}$  of simple roots for $R$ and a set  $\{\omega_i: i\in I\}\subset\lie h^*$    of fundamental weights, i.e., $(\omega_i,d_j\alpha_j)=\delta_{i,j}$, where $\delta_{i,j}$ is the Kronecker delta symbol.
Let $Q$ (resp. $Q^+$) be  the integer span (resp. the nonnegative integer span) of $\{\alpha_i: i\in I\}$ and similarly define  $P$ (resp. $P^+$) be the $\bz$ (resp. $\bz_+$) span of  $\{\omega_i: i\in I\}$. Set $R^+=R\cap Q^+$ and for $\alpha=\sum_{j=1}^nm_j\alpha_j\in R^+$, let $\Ht\alpha=\sum_{j=1}^nm_j$. Denote by  $\theta\in R^+$ be the highest root in $R$ and recall that $[x^\pm_\theta,\lie n^\pm]=0$.
For $i\in I$, we  write $x^\pm_i$, $h_i$, $d_i$ for $x^\pm_{\alpha_i}$, $h_{\alpha_i}$, $d_{\alpha_i}$.

\vskip 6pt

\noindent For $\alpha\in R$, let $\lie g_\alpha$ be the corresponding root space of $\lie g$ and  set $\lie n^\pm=\bigoplus_{\alpha\in R^+}\lie g_{\pm\alpha}$.
Fix non--zero elements  $x^\pm_\alpha\in\lie g_{\pm\alpha}$,  such that $$[h_\alpha, x^\pm_\alpha]=\pm 2x^\pm_\alpha,\ \qquad \  [x^+_\alpha, x^-_\alpha]=h_\alpha,$$ and denote the corresponding subalgebra as $\lie{sl}_2(\alpha)$.

\noindent The next  result is elementary, but we  record it  since it is  important for  Sections 3 and 4.
\begin{lem}\label{rootco} Suppose that $ \beta,\gamma \in R$ are  such that $\beta+\gamma\in R$. Then $\beta+\gamma$ is long if $\beta$ and $\gamma$ are long.\hfill\qedsymbol
\end{lem}

\subsection{} We shall also need, \begin{lem}\label{dalpha} Suppose that   $\lambda\in P^+$ is such that $\lambda=\sum_{i\in I}d_is_i\omega_i$ for some $s_i\in \bz_+$. Then for $\alpha\in R^+$, there exists $s_\alpha\in \bz_+$ such that \begin{equation*} \lambda(h_\alpha) =d_\alpha s_\alpha.\end{equation*} Moreover if $\alpha=\beta+\gamma$  for some $\beta,\gamma\in R$ then $s_\alpha=s_\beta+s_\gamma$. \end{lem}
\begin{pf}
Proceed by induction on $\Ht\alpha$. If $\Ht\alpha=1$, then $\alpha=\alpha_i$ for some $i\in I$ and  $\lambda(h_i)= s_id_i$ showing that induction begins. If $\Ht\alpha>1$,  we can write $\alpha=\beta+\gamma$ with $\beta,\gamma\in R^+$ so that the inductive hypothesis applies to $\beta$ and $\gamma$. If all three roots have the same length, then $h_\alpha=h_\beta+h_\gamma$ and the inductive step is immediate. If $\alpha$ is long then $d_\alpha=1$ and there is nothing to prove.
 If $\alpha$ is short, then by Lemma \ref{rootco} we may assume without loss of generality that $\beta$ is short and $\gamma$ long. This gives,
 $$h_\alpha= h_\beta+d_\alpha h_\gamma.$$ It follows from the inductive hypothesis applied to $\beta$ that $\lambda(h_\beta)=s_\beta d_\beta$ for some $s_\beta\in\bz_+$ and since $d_\alpha=d_\beta$, we get $\lambda(h_\alpha)=d_\alpha(s_\beta+ s_\gamma)$.\end{pf}

\subsection{} For $\mu\in P^+$, let $V(\mu)$ be the irreducible finite--dimensional $\lie g$--module generated by an element $v_\mu$ with defining relations $$x_i^+v_\mu=0,\ \ h_iv_\mu=\mu(h_i)v_\mu,\ \ (x_i^-)^{\mu(h_i)+1}v_\mu=0,\ \ i\in I.$$ It is well--known that any finite--dimensional $\lie g$--module $V$ is isomorphic to a direct sum of irreducible modules $V(\mu)$, $\mu\in P^+$. Further, we may write, $$V=\bigoplus_{\nu\in P} V_\nu,\ \ \ \ V_\nu=\{v\in V: hv=\nu(h)v,\ \ h\in\lie h\},$$ and we set $\wt V=\{\nu\in P: V_\nu\ne 0\}$.

\subsection{} We conclude this section with  the following Lemma which is needed in later sections.
\begin{lem}\label{forcormax2} Suppose that $V$ is a $\lie g[t]$--module and $v\in V$ is such that $$(x_i^-\otimes t^{s_i})v=0,$$ for all $i\in I$ and some $s_i\in\bz_+$. Set $\lambda=\sum_{i\in I}d_is_i\omega_i$. For all $\alpha\in R^+$, we have, $$(x_\alpha^-\otimes t^{s_\alpha})v=0,\ \ s_\alpha\in\bz_+\ \ {\rm{such \ \ that }}\ \  \lambda(h_\alpha)=d_\alpha s_\alpha.$$

\end{lem}
\begin{pf} We prove this by induction on $\Ht\alpha$, with induction beginning for the simple roots by assumption. For $\Ht\alpha>1$, choose $\beta,\gamma\in R^+$ such that $\alpha=\beta+\gamma$, in which case using  Lemma \ref{dalpha}  we see that $s_\alpha=s_\beta+s_\gamma$. Since  $$x_\alpha^-\otimes t^{s_\alpha}= a[x^-_\beta\otimes t^{s_\beta}\ ,\  x_\gamma^-\otimes t^{s_\gamma}]$$ for some non--zero complex number $a$, the inductive step follows and proves the Lemma.
\end{pf}

\section{The modules $V((\xi^\alpha)_{\alpha\in R^+})$} In this section, we first recall the definition of a local Weyl module $W_{\loc}(\lambda)$, $\lambda\in P^+$ and introduce  non--zero graded quotient of these modules which are indexed by an $|R^+|$--tuple of partitions $\bxi=(\xi^\alpha)_{\alpha\in R^+}$ satisfying  $|\xi^\alpha|=\lambda(h_\alpha)$ for all $\alpha\in R^+$. We then give two further equivalent presentations of these modules.

\subsection{}\label{localweylxi}  The definition of the local Weyl modules was given originally in \cite{CPweyl} and later in \cite{CFK} and \cite{FL}. For $\lambda\in P^+,$ the local Weyl module $W_{\loc}(\lambda)$, is the $\lie g[t]$--module generated by an element $w_\lambda$ with defining relations:  for $i\in I$ and $s\in\bz_+$,  \begin{gather}\label{locweylg}(x^+_i\otimes\bc[t])w_\lambda=0,\ \  (h_i\otimes t^s)w_\lambda=\lambda(h_i)\delta_{s,0}w_{\lambda},\\ \label{integloc} (x_i^-\otimes 1)^{\lambda(h_i)+1}w_\lambda=0.\end{gather}  It is simple to see that $\wt W_{\loc}(\lambda)\subset\lambda-Q^+$ and that $\dim W_{\loc}(\lambda)_\lambda=1$. It follows that  $W_{\loc}(\lambda)$ is an indecomposable module and {{that $W_{\loc}(0)$ is isomorphic to  the trivial $\lie g[t]$--module}}.
It was proved in \cite{CPweyl} (see also \cite{CFK}) that the local Weyl modules are  finite--dimensional and so, in particular, we have \begin{equation}\label{intega}(x^-_\alpha\otimes 1)^{\lambda(h_\alpha)+1}w_\lambda=0,\ \ \ \alpha\in R^+.\end{equation} The local Weyl module is clearly graded by $\bz_+$, once we declare the grade of $w_\lambda$ to be zero, and $$W_{\loc}(\lambda)[0]\cong_{\lie g} V(\lambda).$$
Moreover,  $\ev_0 V(\lambda)$ is the unique graded irreducible quotient of $W_{\loc}(\lambda)$.
\vskip 6pt

\subsection{} \label{gar2a} {{ Given  a non--zero element $\lambda\in P^+$, we say that
 $\bxi=(\xi^{\alpha})_{\alpha\in R^+}$ is a $\lambda$--compatible $|R^+|$--tuple of partitions, if  $$\xi^{\alpha}=(\xi^{\alpha}_{\scriptsize 1}\ge\cdots\xi^{\alpha}_s\ge\cdots\ge 0),\ \ |\xi^\alpha|=\sum_{j\ge 1}\xi^{\alpha}_j=\lambda(h_\alpha). $$
Define
 $V(\bxi)$ to be the graded quotient of $W_{\loc}(\lambda)$ by the  submodule generated by the graded elements
 \vskip6pt
  \begin{equation*}
\{ (x_{\alpha}^+\otimes t)^s(x_{\alpha}^-\otimes 1)^{s+r}w_\lambda:\alpha\in R^+, \ s,r\in\bn,  \    s+r\ge 1+ rk+\sum_{j\ge k+1}\xi^{\alpha}_j,\ \ {\rm{for\ some}}\ k\in\bn\}.\end{equation*} \vskip 6pt

   \noindent Denoting by   $v_\bxi$  the image of $w_\lambda$ in $V(\bxi)$, it is clear that
$V(\bxi)$ is  the $\lie g[t]$--module generated by $v_\bxi$ with  defining relations :
\begin{gather}\label{vxi}\lie n^+[t]v_\bxi=0,\ \  (h\otimes t^s)v_\bxi=\delta_{s,0}\lambda(h) v_\bxi,\ \ h\in\lie h,\ \ s\in\bz_+\\ \label{integ} \ (x_{i}^-\otimes 1)^{\lambda(h_i)+1}v_\bxi=0,\ \ i\in I\\ \label{gar2}
(x_{\alpha}^+\otimes t)^s(x_{\alpha}^-\otimes 1)^{s+r}v_\bxi=0,\ \alpha\in R^+, \ s,r\in\bn,  \ \
 s+r\ge 1+ rk+\sum_{j\ge k+1}\xi^{\alpha}_j,\ \ {\rm{for\ some}}\ \  k\in\bn.\end{gather}}}

\vskip 6pt

{\em {For the rest of this section we shall work with a fixed non--zero element  $\lambda\in P^+$ and a fixed  $\lambda$--compatible $\bxi=(\xi^\alpha)_{\alpha\in R^+}$ tuple of partitions.}}

\subsection{} To prove that $V(\bxi)$ is non--zero and to  give alternate formulations of \eqref{gar2} we need some more information.  For $s,r\in\bz_+$, let \begin{equation}\label{maxs}\bs(r,s)=\left\{(b_p)_{p\ge 0}: b_p\in\bz_+, \ \  \sum_{p\ge 0} b_p=r,\ \ \sum_{p\ge 0} pb_p=s\right\}. \end{equation}
 Notice that  $\bs(0,s)$ is the empty set if $s>0$ and  that $$(b_p)_{p\ge 0}\in \bs(r,s)\implies b_p=0\ \ {\rm{if}}\ \ p>s.$$ In particular, $\bs(r,s)$ is finite.
 \vskip6pt

 Given  $x\in\lie g$ and   $s,r\in\bz_+$,
define  elements $\bx(r,s)\in\bu(\lie g[t])$ by, \begin{equation}\label{xlm}\bx(r,s)=\sum_{(b_p)_{p\ge 0}\in\bs(r,s)}(x\otimes 1)^{(b_0)}(x\otimes t)^{(b_1)}\cdots (x\otimes t^s)^{(b_s)},\end{equation} where for any integer $p$ and any $x\in\lie g[t]$,  we set $x^{(p)}=x^p/p!$. Here we understand that $\bx(r,s)=0$ if $\bs(r,s)$ is the empty set.
In particular,  \begin{equation}\label{special1} \bx(0,s)=\delta_{s,0},\qquad \bx(1,s)=x\otimes t^s.\end{equation}
The following result was proved  in \cite{G} (see \cite[Lemma 1.3]{CPweyl} for the current formulation).

\begin{lem}\label{gar} {{Given  $s\in\bn $}}, $r\in\bz_+$ and $\alpha\in R^+$ we have, $$(x_{\alpha}^+\otimes t)^{(s)}(x_{\alpha}^-\otimes 1)^{(s+r)}-(-1)^{s}\bx^-_{\alpha}(r,s)\in\bu(\lie g[t])\lie n^+[t]\bigoplus\bu(\lie n^-[t]\oplus\lie h[t]_+)\lie h[t]_+.$$\hfill\qedsymbol
\end{lem}

Using the Lemma we see immediately that,$$\left((x_{\alpha}^+\otimes t)^{(s)}(x_{\alpha}^-\otimes 1)^{(s+r)} -(-1)^s\bx^-_{\alpha}(r,s)\right) v_\bxi=0, $$
and hence  \eqref{gar2} is equivalent to: for ${\alpha}\in R^+$, $s,r,k\in\bn$,
\begin{equation}\label{gar3} \bx^-_{\alpha}(r,s)v_\bxi =0,\ \  {\rm{if}}\ s+r\ge 1+rk+\sum_{j\ge k+1}\xi^{\alpha}_j.\end{equation} From now on,  we shall  use  freely both presentations of $V(\bxi)$.

\subsection{}{{ Given $\lambda\in P^+$, let $\{\lambda\}=\{(\lambda(h_\alpha)_{\alpha\in R^+}\}$ be the $|R^+|$--tuple of partitions where each partition has at most one part.  We  now prove, \begin{prop}\label{ev0}
\begin{enumerit}
\item[(i)]
 The module  $\ev_0V(\lambda)$ is the unique irreducible  quotient of $V(\bxi)$ and hence  $V(\bxi)$ is a non--zero, indecomposable $\lie g[t]$--module.
 \item[(ii)] We have an isomorphism, $$\ev_0V(\lambda)\cong_{\lie g[t]} V(\{\lambda\}).$$  \end{enumerit}\end{prop}
 \begin{pf} Part (i)  follows if we prove that $v_\lambda$ satisfies \eqref{gar3}. If $s,r \in\bn$,  then $$(b_p)_{p\ge 0}\in\bs(r,s)\implies b_p>0\ \ {\rm{for\ some}}\ \ p>0, $$ and hence
  $$\bx^-_\alpha(r,s)\in\bu(\lie n^-[t])\lie n^-[t]_+,\ \ s,r\in\bn.$$ Since $\lie g[t]_+\ev_0V(\lambda)=0$ we get that $\bx^-_\alpha(r,s)v_\lambda=0$ for $s,r\in\bn$ as required. To prove part (ii)  recall that  $\theta\in R^+$ is the highest root and notice that   the    relation $$\bx_\theta^-(1,1)v_{\{\lambda\}}=0= (x^-_\theta\otimes t)v_{\{\lambda\}},$$ holds in $V(\{\lambda\})$ by taking $s=r=k=1$ in the defining relations. Since $[x^-_\theta,\lie n^-]=0$ and $V(\{\lambda\})=\bu(\lie n^-[t])v_{\{\lambda\}}$ it now follows that  $$(x^-_\theta\otimes t)V(\{\lambda\})=0.$$  Applying elements of $\lie n^+[t]\oplus\lie h[t]_+$ repeatedly, it is now straightforward to see that $$\lie g[t]_+V(\{\lambda\})=0,\ \ i.e.,\ \  \ev_0V(\lambda)\cong_{\lie g[t]} V(\{\lambda\}). $$
 \end{pf}}}

\subsection{} The third presentation requires  an alternate description of the set $\bs(r,s)$. For $k\in\bz_+$, let $\bs(r,s)_k$ (resp. $_k\bs(r,s)$) be the subset of $\bs(r,s)$ consisting of elements $(b_p)_{p\ge 0}$, satisfying $$b_p=0,\ \ p\ge k,\ \ ({\rm{resp.}}\ \ b_p=0\ \ p<k).$$ Let $\ell,m\in\bz_+$ be such that $\bs(r-\ell,s-m)_k\ne\emptyset$ and $ {}_{k}\bs(\ell,m)\ne\emptyset$. Then, the  function \begin{gather*}\bs(r-\ell,s-m)_k\ \times\  {}_{k}\bs(\ell,m)\longrightarrow \bs(r,s), \\ ((b_p)_{p\ge 0}\ ,\ (c_p)_{p\ge 0})\longrightarrow (b_0,\cdots, b_{k-1},c_k,\cdots ),\end{gather*} is clearly one--one and by abuse of notation, we denote the image of this map by  $\bs(r-\ell,s-m)_k\times {}_{k}\bs(\ell,m)$. Clearly,
\begin{equation}\label{alts}\bs(r,s)={} _k\bs(r,s)\bigsqcup \left(\bs(r-r',s-s')_k\times {}_{k}\bs(r',s')\right), \end{equation} where the union is over all pairs $r',s'\in\bz_+$ with $\bs(r-r',s-s')_k\ne\emptyset$ and ${}_k\bs(r',s')\ne\emptyset$.

For $x\in\lie g$, define elements  $\bx(\ell,m)_k$ and $_k\bx(\ell,m)$ of $\bu(\lie g[t])$, \begin{gather*}\bx(\ell,m)_k=\sum_{(b_p)_{\ge 0}\in\bs(\ell,m)_k}(x\otimes 1)^{(b_0)}\cdots (x\otimes t^{k-1})^{(b_{k-1})},\\ _{k}\bx(\ell,m)=\sum_{(b_p)_{\ge 0}\in _k\bs(\ell,m)}(x\otimes t^k)^{(b_k)}\cdots (x\otimes t^m)^{(b_m)}.\end{gather*}
As always, we understand that the element $\bx(\ell,m)_k$ etc. is zero if $S(\ell,m)_k$ etc.  is the empty set. In particular, it follows that\begin{equation}\label{specialcase} \bx(\ell,m)_k\ne 0\implies m\le (k-1)\ell,\qquad _{k}\bx(\ell,k\ell)=(x\otimes t^k)^{(\ell)} .\end{equation}

  \begin{lem}\label{rearrange} Suppose that $s,r,k\in\bn$ are such that $s+r\ge kr+K$ for some $K\in\bz_+$.
  Then,\begin{equation*}
  \bx(r,s)={} _{k}\bx(r,s)+ \sum\bx(r-r',s-s')_k\  \ {}_k\bx(r',s'),\end{equation*} where the sum is over all pairs $r',s'\in\bz_+$ satisfying $r'<r$, $s'\le s$ {\em and} $s'+r'\ge r'k+K$.
  \end{lem}
  \begin{pf} It is trivial to see from the alternative description of the set $\bs(r,s)$ given  in \eqref{alts}, that \begin{equation*}
  \bx(r,s)={} _k\bx(r,s)+ \sum\bx(r-r',s-s')_k\  \ {}_k\bx(r',s'),\end{equation*} where the sum is over all $r'<r$ and $s'\le s$. Using \eqref{specialcase} we see that we can also assume that $s-s'\le (k-1)(r-r')$. Since $$s+r=(s-s')+(r-r')+s'+r'\ge k(r-r')+kr'+K,$$ we get $s'+r'\ge kr'+K$ as required.
  \end{pf}

\subsection{} We now  prove, \begin{prop}\label{max}  Let $V$ be any representation of $\lie g[t]$ and let $v\in V$, $x\in\lie g$ and $K\in\bz_+$. Then,
\begin{gather*}\bx(r,s)v=0\ {\rm{for\ all}}\ s,r,k\in\bn
 \ \ {\rm{with}}\  s+r\ge 1+kr+K\ \iff\\  _{k}\bx(r,s)v=0  \ {\rm{ for \ all}}\  s,r,k \in\bn \  {\rm{with}}\ s+r\ge 1+ kr+K.\end{gather*}
\end{prop}
\begin{pf} Suppose first that $\bx(r,s)v=0$ for all $s,r,k\in\bn$ with $s+r\ge 1+kr+K$. We shall prove by induction on $r$ that $_{k}\bx(r,s)v=0$ for all $s,r,k\in\bn$ with $s+r\ge 1+ kr+K$.
 Suppose that $r=1$ and let $s\in\bn$ be such that $s+1\ge 1 +k+K$. Then  we have $s\ge k$ and so  $$\bx(1,s)={} _{k}\bx(1,s)=x\otimes t^s,$$
and hence we have that ${}_k\bx(1,s)v=0$.
 Assume now that we have proved the statement  for all $r'<r$ and all $s'\le s$ with $s'+r'\ge 1+ kr'+K$.
Using Lemma \ref{rearrange}, we get $$0=\bx(r,s)v= {}_{k}\bx(r,s)v  +\left(\sum \bx(r-r',s-s')_k\ {}_k\bx(r',s')\right)v,$$  where the sum is over all $r'<r$ and $s'\le s$ with $r'+s'\ge 1+kr'+K$. The inductive hypothesis applies to the second term on the right hand side and hence we get  ${}_{k}\bx(r,s)v=0$.
The converse statement is obvious by using Lemma \ref{rearrange}.
\end{pf}
\subsection{} We can now give the third presentation of $V(\bxi)$.
\begin{prop}\label{third} The module $V(\bxi)$ is generated by the element $v_\bxi$ with defining relations \eqref{vxi},\eqref{integ} and,
\begin{equation}\label{third2} _{k}\bx^-_{\alpha}(r,s)v_\bxi=0,\ \ {\alpha}\in R^+,\ \ s,r,k\in\bn, \ \ s+r\ge 1+kr+\sum_{j\ge k+1}\xi^{\alpha}_j.\end{equation}
 In particular, for all ${\alpha}\in R^+$, $r, k\in\bn$with  $r\ge 1+ \sum_{j\ge k+1}\xi^{\alpha}_j$, we have \begin{equation}\label{demreltype}(x^-_{\alpha}\otimes t^k)^{r}v_\bxi=0.\end{equation}
\end{prop}
\begin{pf} The first statement is immediate from  Proposition \ref{max}. For the second, take  $s=kr$. Then we have $s+r\ge 1+kr+ \sum_{j\ge k+1}\xi^{\alpha}_j$ and  using the  second equation in  \eqref{specialcase}, we get $$ _{k}\bx^-_\alpha(r, kr)v_\bxi= (x^-_\alpha\otimes t^k)^{r}v_\bxi=0.$$\end{pf}

\begin{cor} If $|\xi^\alpha|>0$ and $s_{\alpha}$ is the number of parts of $\xi^{\alpha}$, then $$(x_{\alpha}^-\otimes t^p)v_\bxi=0,\ p\ge s_{\alpha}. $$\hfill\qedsymbol
\end{cor}

\section{The connection with  Demazure modules. }

In this section, we study very special kinds of $\lambda$--compatible partitions and prove that in this case,
   the defining relations of $V(\bxi)$ can be greatly simplified. This allows us to make  connections with  well--known $\lie g[t]$--modules, such as the local Weyl modules, Kirillov--Reshetikhin modules and  the  Demazure modules of level $\ell$. Our results are new when $\ell>1$ and are  much simpler than the original proofs  (see  \cite{CL},\cite{CPweyl},\cite{FoL},\cite{Naoi}) even when $\ell=1$. We shall use freely and without comment, the notation established in the earlier sections.

\subsection{}  It will be convenient to use another standard notation for partitions.  Namely, if  $i_1>\cdots> i_r$ are the  distinct  non--zero parts of the partition and  $i_k$ occurs $s_k$ times
 then  we shall denote this partition by $(i_1^{s_1},\cdots, i_r^{s_r} )$.
A  partition $\xi$ is said to be rectangular if it is either the empty partition or of the form $(k^m)$ for some $k,m\in\bn$. A partition is said to be a fat hook if it is of the form $(k_1^{s_1}, k_2^{s_2})$ with $k_j,s_j\in\bn$, $j=1,2$. We shall say that the fat hook is special if $s_2=1$.

\begin{thm}\label{genmax2} Let $\bxi=(\xi^{\alpha})_{\alpha\in R^+}$ be a $\lambda$--compatible $|R^+|$--tuple of partitions. Assume  that  $\xi^\alpha $ is  either rectangular or a special
 fat hook for $\alpha\in R^+$. Then, $V(\bxi)$ is isomorphic to the quotient of $W_{\loc}(\lambda)$ by the submodule generated by the elements,

 \begin{gather}\label{finer1}\{ (x^-_{\alpha}\otimes t^{s_{\alpha}})w_\lambda:\alpha\in R^+\}\ \bigcup\   \{(x^-_{\alpha}\otimes t^{s_{\alpha}-1})^{\xi^\alpha_{s_\alpha}+1}w_\lambda:\alpha\in R^+, \ \ \xi^{\alpha}\ \ {\rm{{a\ special\ fat \ hook}}}\},
\end{gather} \vskip 6pt
\noindent where $s_\alpha$ is zero if $\xi^\alpha$ is the empty partition and  is   the number of non--zero parts of $\xi^\alpha$ otherwise.

\end{thm}

\begin{pf} Let $U$ be the submodule of $W_{\loc}(\lambda)$ generated by the elements in \eqref{finer1} and
let $\tilde V(\bxi)$ be the corresponding quotient of $W_{\loc}(\lambda)$.  Using Proposition \ref{third} and taking   $k\in\{s_\alpha-1,s_\alpha\}$ in equation \eqref{demreltype}, we see that
   that $V(\bxi)$ is a quotient of $\tilde V(\bxi)$. To prove that they are isomorphic, we must  show  that:  for    $\alpha\in R^+$  and $k,s,r\in\bn$, either \begin{equation}\label{claim1}  s+r\ge 1+ rk+\sum_{j\ge k+1}\xi^{\alpha}_j\implies
  (x_{\alpha}^+\otimes t)^s(x_{\alpha}^-\otimes 1)^{s+r}w_\lambda\in U,\end{equation} or \begin{equation}\label{claim2}  s+r\ge 1+ rk+\sum_{j\ge k+1}\xi^{\alpha}_j\implies
  \bx^-_\alpha(r,s)w_\lambda\in U. \end{equation}

\vskip 6pt

 If $r\ge \xi^{\alpha}_{\scriptsize{1}}$, then $$s+r\ge 1+k\xi_1^{\alpha}+\sum_{j\ge k+1}\xi^{\alpha}_j\ge 1+\sum_{j\ge 1}\xi^{\alpha}_j=|\xi^{\alpha}|+1. $$ By \eqref{intega}, we know $(x_{\alpha}^-\otimes 1)^{s+r}w_\lambda=0$ and so equation \ref{claim1} is proved in this case.
 \vskip 6pt

 \noindent If  $r<\xi^{\alpha}_{\scriptsize{1}}$  we shall prove that \eqref{claim2} is satisfied. Observe that if $(b_p)_{p\ge 0}\in\bs(r,s)$ is  such that $b_m>0$ for some $m\ge s_\alpha$, then   $$\left((x_\alpha^-\otimes 1)^{(b_0)}\cdots (x_\alpha^-\otimes t^{m})^{(b_m)}\cdots\right)w_\lambda\in U,$$ and hence, we get
 \begin{equation}\label{crux}\left( \bx^-_\alpha(r,s)-\bx^-_\alpha(r,s)_{s_\alpha}\right)w_\lambda\in U.\end{equation}

If $r\le \xi^\alpha_{s_\alpha}$, we claim that,   \begin{equation}\label{inequality} s+r\ge 1+kr+\sum_{j\ge k+1}\xi^{\alpha}_j\implies s+r\ge 1+s_{\alpha}r.\end{equation}
 For the claim, notice there is nothing to prove if $k\ge s_{\alpha}$, and    if $k< s_{\alpha}$, then $$s+r\ge 1+kr+ \sum_{j\ge k+1}\xi^{\alpha}_j\ge 1+kr+(s_{\alpha}-k)\xi_{s_\alpha}^{\alpha}\ge 1+kr+(s_{\alpha}-k)r.$$
 This means that if    $(b_p)_{p\ge 0}\in\bs(r,s)$, then we must have $b_m>0$ for some $m\ge s_{\alpha}$,
since otherwise we would have $s=\sum_{p<s_{\alpha}}pb_p\le r(s_{\alpha}-1)$. In particular, we get $\bx^-_\alpha(r,s)_{s_\alpha}=0$ and  equation \eqref{crux} now proves that $\bx^-_\alpha(r,s)w_\lambda\in U$. This  also completes the proof when $\xi^\alpha$ is rectangular.

 \vskip6pt

 If  $\xi^\alpha$ is a special fat hook  we have to also consider the case when
 $\xi_1^{\alpha}>r> \xi^\alpha_{s_\alpha}$.
If $(b_p)_{p\ge 0}\in\bs(r,s)_{s_{\alpha}}$ then $s=(s_{\alpha}-1)b_{s_{\alpha}-1}+\cdots +b_1$.  If we  prove that  \begin{equation}\label{crux2} b_{s_{\alpha}-1}\ge \xi^{\alpha}_{s_\alpha}+1,\end{equation} then it obviously follows that $\bx^-_\alpha(r,s)_{s_\alpha} w_\lambda\in U$. Using \eqref{crux} again, we will have $\bx_\alpha ^-(r,s) w_\lambda\in U$ which would complete the proof of the theorem.

  To prove equation \ref{crux2}, observe that $$(s_{\alpha}-1)b_{s_{\alpha}-1}+(r-b_{s_{\alpha}-1})(s_{\alpha}-2)\ge s\ge 1+ r(k-1)+\sum_{j\ge k+1}\xi^{\alpha}_j,$$ and hence
$$\ \ b_{s_{\alpha}-1}\ge 1+ r(k-s_{\alpha}+1)+\sum_{j\ge k+1}\xi^{\alpha}_j.$$
Since $r>\xi^\alpha_{s_\alpha}$, we see that equation \ref{crux2} is immediate if  $k\ge s_{\alpha}$. If $k< s_{\alpha}$, then
$$b_{s_{\alpha}-1}\ge 1+\sum_{s_{\alpha}> j\ge k+1}(\xi^{\alpha}_j-r)+\xi^\alpha_{s_\alpha}\ge 1+ \xi^\alpha_{s_\alpha},$$ where the last inequality is because $\xi^\alpha$ is a special fat hook, which means $\xi^\alpha_1 =\xi^\alpha_j>r$ for all $1\le j\le s_\alpha-1$.

\end{pf}

 \subsection{} Given $\ell\in\bn$ and  a non--zero element $\lambda\in P^+$, we now define in a canonical way, an $|R^+|$--tuple  of  $\lambda$--compatible partitions.  For $\alpha\in R^+$, with $\lambda(h_\alpha)>0$,  let  $s_\alpha,m_\alpha\in\bn$ be the unique positive integers so that  $$\lambda(h_{\alpha})=(s_{\alpha}-1)d_{\alpha}\ell+m_{\alpha},\ \ \ 0<m_{\alpha}\le d_{\alpha}\ell,$$ where we recall that $d_\alpha=2/(\alpha,\alpha)$. If $\lambda(h_\alpha)=0$ set $s_\alpha=0=m_\alpha$.
Let $\bxi(\ell,\lambda)=(\xi^{\alpha})_{\alpha\in R^+}$ be the  $|R^+|$--tuple of partitions given by: $\xi^\alpha $ is the empty partition if $\lambda(h_\alpha)=0$ and otherwise, is the  partition $((d_\alpha\ell)^{s_\alpha-1}, m_\alpha)$.  In particular, $\xi^\alpha$ is either rectangular or a special fat hook. In the rest of this section we use Theorem \ref{genmax2} to analyze the modules $V(\bxi(\ell,\lambda))$.

 \subsection{}  The first result is,
\begin{lem}\label{first} Let $\ell\in\bz_+$ and suppose that $\lambda=\ell\left(\sum_{i\in I}d_is_i\omega_i\right)$. Then, $V(\bxi(\ell,\lambda))$ is the quotient of $W_{\loc}(\lambda)$  by the submodule generated by the elements $$\{(x_i^-\otimes t^{s_i})w_\lambda:\  i\in I,\ \ s_i>0\}.$$
\end{lem}
\begin{pf}  Applying  Lemma \ref{dalpha} to $(1/\ell)\lambda$ proves that $\lambda(h_\alpha)= d_\alpha\ell s_\alpha$ for some $s_\alpha\in R^+$, i.e, that $\xi^\alpha$ is rectangular for all $\alpha\in R^+$. Hence by Theorem \ref{genmax2}, it suffices to prove that
\begin{equation}\label{1}(x_i^-\otimes t^{s_i})v_{\bxi(\ell,\lambda)}=0, \ i\in I\implies (x_\alpha^-\otimes t^{s_\alpha})v_{\bxi(\ell,\lambda)}=0,\ \alpha\in R^+.\end{equation} Recalling that  $(x_i^-\otimes 1)w_\lambda=0$ if $s_i=0$ equation \eqref{1} follows by applying Lemma \ref{forcormax2}  to $(1/\ell)\lambda$.
\end{pf}
\subsection{}  Our next result considers the case when $\lambda\in P^+$ and $\ell=1$.
\begin{prop}\label{second} For  $\lambda\in P^+$, the module $V(\bxi(1,\lambda))$ is the quotient of the local Weyl module  by the submodule generated by the elements in the following two sets:  \begin{gather}\label{three} \{(x^-_\alpha\otimes t^{s_\alpha})w_\lambda: \alpha\in R^+\ {\rm{with }}\  d_\alpha>1\},\\
\label{four}\{ (x^-_\alpha\otimes t^{s_\alpha-1})^{2}w_\lambda:\alpha\in R^+ \ {\rm{with }}\ d_\alpha=3\ \ {\rm{and}}\ \  m_\alpha=1\}.\end{gather}
  In particular, if $\lie g$ is simply laced, then $W_{\loc}(\lambda)\cong V(\bxi(1,\lambda))$.
\end{prop}

\begin{pf} Let $U$ be the submodule generated by the elements in \eqref{three} and  \eqref{four}. If $\lambda(h_\alpha)=0$, then  $(x_{\alpha}^-\otimes 1)w_\lambda=0$ and there is nothing to prove. Assume now that $\lambda(h_\alpha)>0$.  Since  $\ell=1$, we have $0<m_\alpha\le d_\alpha$ and $\lambda(h_\alpha)=d_\alpha(s_\alpha-1)+m_\alpha$. In particular if $d_\alpha=1$ then $\xi^\alpha$ is rectangular and $s_\alpha=\lambda(h_\alpha)$. Applying Theorem \ref{genmax2} to $V(\bxi(1,\lambda))$, we see that the result follows if we prove that: \begin{gather*}d_\alpha=1\implies (x_\alpha^-\otimes t^{s_\alpha})w_\lambda=0, \\ d_\alpha> 1\ {\rm{and}}\ \  m_\alpha=d_\alpha-1\implies (x_\alpha^-\otimes t^{s_\alpha-1})^{d_\alpha}w_\lambda\in U.\end{gather*}
The first statement follows by using \eqref{intega} and Lemma \ref{gar} to see that $$0=(-1)^{s_{\alpha}}(x_\alpha^+\otimes t)^{(s_\alpha)}(x_\alpha^-\otimes 1)^{(s_\alpha+1)} w_\lambda =(x_\alpha^-\otimes t^{s_\alpha})w_\lambda=\bx^-_\alpha(1,s_\alpha)w_\lambda.$$

 Suppose that $d_\alpha>1$ and $m_\alpha=d_\alpha-1$, in which case $\lambda(h_\alpha)=s_\alpha d_\alpha-1$ and we get by using \eqref{intega} and Lemma \ref{gar} again, that  \begin{gather*}  0=(-1)^{d_{\alpha}(s_{\alpha}-1)}(x^-_\alpha\otimes t)^{(d_\alpha(s_\alpha-1))}(x_\alpha^-\otimes 1)^{(d_\alpha s_\alpha)}w_\lambda\\ = \bx^-_\alpha(d_\alpha, d_\alpha (s_\alpha-1))w_\lambda
=  (x_\alpha^-\otimes t^{s_\alpha-1})^{d_\alpha}w_\lambda + Xw_\lambda,\end{gather*} where $X$ is in the left ideal of $\bu(\lie n^-[t])$ generated by the elements $(x^-_\alpha\otimes t^p)$ with $p\ge s_\alpha$. Hence $Xw_\lambda\in U$ which give $$(x_\alpha^-\otimes t^{s_\alpha-1})^{d_\alpha}w_\lambda\in U.$$
\end{pf}

   \subsection{} We now  recall from \cite{FoL} and \cite{Naoi} the definition of the Demazure modules which are appropriate to our study.
 For $\ell\in\bz_+$ and  $\lambda\in P^+$  the Demazure module $D(\ell,\lambda)$  is the graded quotient of $W_{\loc}(\lambda)$ generated by the elements,
 \begin{gather} \label{demrelo}
 \{(x^-_\alpha\otimes t^p)^{r+1}w_\lambda:\ \ p\in\bz_+,\ \   r\ge\max\{0,\lambda(h_\alpha)-d_\alpha\ell p\},\ \ \text{for} \ \alpha\in R^+\}. \end{gather}

 \begin{thm}\label{genreldem1} Let   $\ell\in\bz_+$ and $\lambda\in P^+$.
  We have an isomorphism of graded  $\lie g[t]$--modules   $V(\bxi(\ell,\lambda)) \cong D(\ell,\lambda)$.
   Equivalently, $D(\ell,\lambda)$ is the quotient of $W_{\loc}(\lambda)$ by the submodule generated by the elements
 \begin{gather}\label{demrel}\{ (x^-_{\alpha}\otimes t^{s_{\alpha}})w_\lambda: \alpha\in R^+\}\ \bigcup\  \{ (x^-_{\alpha}\otimes t^{s_{\alpha}-1})^{m_{\alpha}+1}w_\lambda: \alpha\in R^+, \ m_\alpha<d_\alpha\ell\}\ .\end{gather}\end{thm}
 \begin{pf} Let $U'$ be the submodule of $W_{\loc}(\lambda)$ generated by the elements in \eqref{demrelo}. Taking $p\in\{s_\alpha, s_\alpha-1\}$ in \eqref{demrelo} and recalling that $\lambda(h_\alpha)=d_\alpha\ell(s_\alpha-1)+m_\alpha$, we see   that the elements $(x_\alpha^-\otimes t^{s_\alpha})w_\lambda$ and $(x^-_\alpha\otimes t^{s_\alpha-1})^{m_\alpha+1}$ are in $U'$. Theorem \ref{genmax2} implies that  the canonical map $W_{\loc}(\lambda)\to D(\ell,\lambda)$ factors through to a map
 of $\lie g[t]$--modules $$V(\bxi(\ell,\lambda))\to D(\ell,\lambda).$$ To prove that it is an isomorphism we must prove that the additional defining relations of $D(\ell,\lambda)$ hold in $V(\bxi(\ell,\lambda))$.
Namely that for all $\alpha\in R^+$  and $p\in\bz_+$, we have
\begin{equation}\label{genrelationmainthm}
(x^-_\alpha\otimes t^p)^{r+1}v_{\bxi(\ell,\lambda)}=0,\ \ \  r\ge\max\{0,\ \lambda(h_\alpha)-d_{\alpha}\ell p\}.
\end{equation}  If    $p\ge s_{\alpha}$, we have   $r\ge 0$ while, if $p< s_{\alpha}$, then $$r\ge\lambda(h_\alpha)-d_{\alpha}\ell p=d_{\alpha}\ell (s_{\alpha}-p-1)+m_{\alpha}\ge \sum_{j\ge p+1}\xi^{\alpha}_j.$$ In either case,  \eqref{genrelationmainthm} follows
by taking $k=p$ in \eqref{demreltype} and  the proof of the  Theorem is complete.
\end{pf}

\subsection{} Together with Proposition \ref{second} we have,
 \begin{cor}\label{demloc} The module $D(1,\lambda)$ is the quotient of the local Weyl module  by the submodule generated by the elements in the following two sets:  \begin{gather*}\{(x^-_\alpha\otimes t^{s_\alpha})w_\lambda: \alpha\in R^+\ {\rm{such \ that}}\  d_\alpha>1\}\\
\{ (x^-_\alpha\otimes t^{s_\alpha-1})^{2}w_\lambda:\alpha\in R^+ \ {\rm{such \ that}}\ d_\alpha=3\ \ {\rm{and}}\ \  m_\alpha=1\}.\end{gather*}
  In particular, if $\lie g$ is simply laced, then $W_{\loc}(\lambda)\cong D(1,\lambda)$.\hfill\qedsymbol
 \end{cor}

\subsection{} We end this section with some remarks.  Corollary \ref{demloc} was proved earlier by very different methods.   In \cite[Section 6]{CPweyl} it was shown that $W_{\loc}(\ell\omega)$ is of dimension $2^\ell$ when $\lie g$ is isomorphic to $\lie{sl}_2$.  On the other hand it was known  that $D(1,\ell\omega)$ is also of dimension $2^\ell$ (\cite{FoL2},\cite{KMOTU},\cite{M}) which  shows  that $W_{\loc}(\ell\omega)\cong D(1,\ell\omega)$. This was then used in \cite[Section 3]{CL} and later in \cite[Section 3.4]{FoL} to prove the result for $\ell=1$ and $\lie{sl}_{r+1}$ and for an arbitrary simply laced Lie algebra respectively. The case of the non--simply laced algebras the result  was proved in \cite[Section 4]{Naoi} using a case by case approach.
The methods used in these papers  are more complicated and   do not appear to generalize  to the case of $\ell >1$ and arbitrary $\lambda$.

\section{ Applications to  Fusion Products of Demazure modules}
We begin by recalling   the definition of fusion products of $\lie g[t]$--modules given in \cite{FL}.  We then  prove some elementary results which are  used repeatedly in the rest of the paper. The main result of this section is Proposition \ref{mapsdem} on the fusion product of Demazure modules.

 \subsection{} Suppose that $V$ is a cyclic $\lie g[t]$--module generated by an element $v$. Define a filtration  $F^rV$, $r\in\bz_+$ on $V$, by $$F^rV= \sum_{0\le s\le r}\bu(\lie g[t])[s]v.$$ Clearly each $F^rV$ is a $\lie g$--module and the associated graded space $\gr V$ acquires a natural structure of a cyclic graded $\lie g[t]$--module with action given by $$(x\otimes t^s)(\overline w)=\overline{(x\otimes t^s)w},\ \ \ \overline w\in F^rV/F^{r-1}V.$$ Moreover, $\gr V\cong V$ as $\lie g$--modules and $\gr V$ is the  cyclic $\lie g[t]$--module generated by the image $\bar v$ of $v$ in $\gr V$. The following is trivial but very useful.
 \begin{lem}\label{elemfusion} Let $V$ be a cyclic $\lie g[t]$--module generated by $v\in V $. For all $u\in V$, $x\in\lie g$, $r\in\bn$, $a_1,\cdots, a_r\in\bc$, we have $$(x\otimes t^r)\bar u= (x\otimes (t-a_1)\cdots (t-a_r))\bar u,$$ where $\bar u$ is the image of $u$ in $\gr V$.\hfill\qedsymbol
 \end{lem}

\subsection{} Given any $\lie g[t]$--module $V$ and $z\in\bc$ let  $V^z$ be the $\lie g[t]$--module defined by $$(x\otimes t^r)v=(x\otimes (t+z)^r)v,\ \ x\in \lie g,\ \  r\in\bz_+,\ \ v\in V.$$
 Let  $V_1,\cdots, V_m$ be  finite--dimensional graded $\lie g[t]$--modules generated by elements $v_j$, $1\le j\le m$ and  let  $z_1,\cdots,z_m$ be   distinct complex numbers. Let  $$\bv= V_1^{z_1}\otimes\cdots\otimes V_m^{z_m},$$ be the corresponding product of $\lie g[t]$--modules.
 It is easily checked (see \cite[Proposition 1.4]{FL})  that $$\bv=\bu(\lie g[t])(v_1\otimes \cdots\otimes v_m).$$ The corresponding associated graded $\lie g[t]$--module $\gr \bv$ is called the fusion product of $V_1,\cdots ,V_m$ with parameters $z_1,\cdots, z_m$. It is denoted as $V_1^{z_1}*\cdots *V_m^{z_m}$ and is generated by the image   of $v_1\otimes\cdots\otimes v_m$.

{\em In the rest of the paper, we shall frequently,  for ease of notation,  suppress the dependence of the fusion product on the parameters and just write $V_1*\cdots*V_m$ for $V_1^{z_1}*\cdots *V_m^{z_m}$. But unless explicitly stated, it should be assumed that the fusion product does depend on these parameters}.

\vskip 6pt

Given elements $u_s\in V_s$, $1\le s\le m$, we shall denote by $u_1*\cdots *u_m\in V_1*\cdots*V_m $ the image of the element $u_1\otimes\cdots \otimes u_m\in V_1^{z_1}\otimes\cdots\otimes V_m^{z_m}$.

\subsection{} \begin{lem}\label{fusionweyl} Let $\lambda_s\in P^+$ and  $V_s$ be a $\lie g[t]$--module quotient of $W_{\loc}(\lambda_s)$  for $1\le s\le m$.   Then $V_1*\cdots*V_m$ is a graded $\lie g[t]$--module quotient of $W_{\loc}(\lambda)$, where $\lambda=\sum_{s=1}^m\lambda_s$.\end{lem}
\begin{pf} Let $v_s$ be the image of $w_{\lambda_s}$ in $V_s$. It suffices to prove that  $v_1*\cdots *v_m$ satisfies the defining relations of $w_\lambda$. Note that for $r>0$ and $h\in\lie h$, we have $$(h\otimes t^r)(v_1\otimes \cdots \otimes v_m)= \left(\sum_{s=1}^mz_s^r\lambda_s(h)\right)(v_1\otimes \cdots \otimes  v_m)\in\bu(\lie g)(v_1\otimes\cdots\otimes v_m).$$ Hence if $r>0$, we  see that the image of $(h\otimes t^r)(v_1\otimes \cdots \otimes v_m)$ in the fusion product is zero. The other defining relations of $W_{\loc}(\lambda)$ are trivially satisfied and the proof is complete.
\end{pf}

\subsection{} We shall need  some  results on Demazure modules which were  proved in \cite{FoL2}, \cite{FoL}. We note here, that these papers work with a special family of Demazure modules $D(\ell,\lambda)$. In  the notation of the current paper, they only work with the modules of the form  $D(\ell,\ell\mu)$, where $\mu\in P^+$ is of the form  $\mu=\left(\sum_{i\in I }d_is_i\omega_i\right)\in P^+$.
In view of this, it is convenient to define a subset $\Gamma$ of $\bz_+\times P^+$ as follows: $$\Gamma=\left\{(\ell,\lambda)\in\bz_+\times P^+:\ \   \lambda=\ell\sum_{i\in I}d_is_i\omega_i\right\}.$$ In particular, if $(\ell,\lambda)\in\Gamma$ it follows from Theorem \ref{genreldem1} and Lemma \ref{first} that the Demazure module $D(\ell,\lambda)$,  where $\lambda=\ell\left(\sum_{i\in I}d_is_i\omega_i\right)$ is the quotient of the local Weyl module by the relations $(x_i^-\otimes t^{s_i})w_\lambda$, $i\in I$ with $s_i>0$.
\vskip 6pt

The following was established in \cite[section 4]{FoL}, \cite[section 3]{FoL2}
\begin{thm}\label{demprop} Let  $(\ell,\lambda), (\ell,\mu)\in\Gamma$.\begin{enumerit}
\item[(i)] For all $s\in\bz$, we have  $\dim\Hom_{\lie g[t]}(\tau_s D(\ell, \lambda),D(\ell,\mu))\le 1$ and moreover, any non--zero map is injective.
\item[(ii)] $\dim D(\ell,(\lambda+\mu))=\dim D(\ell,\lambda)\dim D(\ell,\mu)$.\hfill\qedsymbol
\end{enumerit}
\end{thm}

\subsection{} We can now prove,
\begin{prop}\label{mapsdem} Let $(\ell,\lambda)\in\Gamma$ and  suppose  that there exists $(p_j,\mu_j)\in\Gamma$, $1\le j\le m$  such that $$\lambda=\mu_1+\cdots+\mu_m, \qquad \frac1\ell\lambda(h_i)\ge\sum_{j=1}^m\frac1{p_j}\mu_j(h_i),\ \ i\in I.$$ There exists a non--zero surjective map of graded $\lie g[t]$--modules, $$D(\ell, \lambda)\to D(p_1, \mu_1)*\cdots*D(p_m,\mu_m)\to 0,$$
and in   the special case when $p_1=\cdots=p_m=\ell$  we have an isomorphism $$D(\ell, \lambda)\cong D(\ell, \mu_1)*\cdots*D(\ell, \mu_m).$$
 \end{prop}
 The following is immediate.
 \begin{cor} The fusion product of a finite number of modules of the form $D(\ell,\mu)$, $(\ell,\mu)\in\Gamma$ for a fixed $\ell$ is independent of the choice of parameters. \end{cor}
\begin{rem} The second statement of the proposition and the Corollary  was proved earlier in \cite[Section 3.5]{FoL} using a result in \cite{FKL}.
\end{rem}

\noindent{\em{Proof of Proposition \ref{mapsdem}.}} \ Let $v_s\in D(p_s,\mu_s)$ be the image of the generator  $w_{\mu_s}$ of $W_{\loc}(\mu_s)$ for $1\le s\le m$. Using Lemma \ref{fusionweyl}, we see that there exists a surjective map  of graded $\lie g[t]$--modules,
 $$W_{\loc}(\lambda)\to D(p_1, \mu_1)*\cdots*D(p_m,\mu_m)\to 0.$$
Writing $\lambda=\ell\sum_{i\in I}d_is_i\omega_i$, and using  Theorem \ref{genreldem1} and Lemma \ref{first} (see the comments preceding the statement of Theorem \ref{demprop}),
it suffices to show that $$(x_i^-\otimes t^{s_i})(v_1*\cdots*v_m)=0, i\in I, s_i>0.$$ Write $\mu_k=p_k\sum_{i\in I}d_is_{i,k}\omega_i,\ \ 1\le k\le m,$ and note that we are given that $$s_i\ge\sum_{k=1}^m s_{i,k},\ \ i\in I.$$
Setting    $b_i=s_i-\sum_k s_{i,k}$ and taking $z_1,\cdots,z_m$ be the parameters involved in the fusion product, we see that
\begin{gather*}(x^-_i\otimes t^{b_i}(t-z_1)^{s_{i,1}}\cdots(t-z_m)^{s_{i,m}})(v_1\otimes\cdots \otimes v_m) \\
={\sum_{r=1}^m v_1\otimes\cdots\otimes(x^-_i\otimes (t+z_r)^{b_i}(t-z_1+z_r)^{s_{i,1}}\cdots t^{s_{i,r}}\cdots (t+z_r-z_m)^{s_{i,m}})v_r\otimes\cdots\otimes v_m}.\end{gather*} For $1\le r\le m$, we know  that the relation $$(x_i^-\otimes t^b)v_r=0,\ \ b\ge s_{i,r}, $$ holds in $D(p_r,\mu_r)$ and hence we have shown that$$(x^-_i\otimes t^{b_i}(t-z_1)^{s_{i,1}}\cdots(t-z_m)^{s_{i,m}})(v_1\otimes \cdots \otimes v_m)=0. $$It now follows by using Lemma \ref{elemfusion} that $$(x_i^-\otimes t^{s_i})(v_1*\cdots*v_m)=(x_i^-\otimes t^{b_i}(t-z_1)^{s_{i,1}}\cdots(t-z_m)^{s_{i,m}})(v_1*\cdots * v_m)=0,$$ which proves the existence of the map $$D(\ell, \lambda)\to D(p_1, \mu_1)*\cdots*D(p_m,\mu_m)\to 0.$$
The second statement of the proposition  is now immediate by using  Theorem \ref{demprop}(ii).\hfill\qedsymbol

\section{Kirillov--Reshetikhin modules and $Q$--systems}
 In this section, we discuss further consequences of our study and establish the connections with the  $Q$--systems introduced in \cite{HKOTY}.    We will  use freely the notation established  in the previous sections.
\subsection{} We recall from \cite[Section 2]{CMkir} the definition of the Kirillov--Reshetikhin modules. Thus given $i\in I$ and $m\in\bz_+$, the Kirillov--Reshetikhin module $KR(m\omega_i)$ is the quotient of the  $W_{\loc}(m\omega_i)$ by the submodule generated by the element  $(x_i^-\otimes t)w_{m\omega_i}$. We will denote  the image of $w_{m\omega_i}$ in $KR(m\omega_i)$ by $w_{i,m}$. It is trivially checked  that  \begin{equation}\label{bigger}(x_i^-\otimes t^s)w_{i,m}=0,\ \ s\ge 1.\end{equation}

\vskip 6pt

The following is a consequence of Lemma \ref{first}, Theorem \ref{genreldem1} and Proposition \ref{mapsdem}. We remark that the isomorphism between  the Kirillov--Reshetikhin modules and the Demazure modules was proved earlier in \cite[Section 5]{CMkir} and \cite[Section 3.2]{FoL}.
\begin{prop} \label{krdem} For $i\in I$ and $\ell\in\bz_+$, we have an isomorphism of $\lie g[t]$--modules, \begin{gather*}KR(d_i\ell\omega_i)\cong V(\bxi(\ell, d_i\ell\omega_i))\cong D(\ell, d_i\ell\omega_i).\end{gather*} Moreover if $\lambda=\ell(\sum_{i\in I}d_is_i\omega_i)$ then,\begin{gather*}
D(\ell,\lambda)\cong KR(d_1\ell\omega_1)^{*s_1}*\cdots *KR(d_n\ell\omega_n)^{*s_n}\cong V(\bxi(\ell,\lambda)).\end{gather*}\hfill\qedsymbol \end{prop}

\subsection{} We recall the definition of the $Q$--system  given in \cite[Section 7]{HKOTY}. Consider the ring {$\bz[x_1^{\pm 1},\cdots, x_n^{\pm 1}]$  in the indeterminates $x_1,\cdots,x_n$, where we recall that $n$ is the rank of $\lie g$. \iffalse Since $\lie g$ is of classical type,  for any $\lambda\in P^+$, the character of  $V(\lambda)$ can be regarded as  an element of this ring.\fi  A $Q$--system for $\lie g$ is a family of functions $\{Q^{(i)}_m: i\in I, m\in\bz_+\}$  satisfying  $Q^{(i)}_0=1$, and $$Q_m^{(i)}Q_m^{(i)}= Q_{m+1}^{(i)}Q_{m-1}^{(i)}+ \prod_{j\sim i}\prod_{k=0}^{-C_{ij}-1}Q^{(j)}_{\left\lceil\frac{mC_{ji}-k}{C_{ij}}\right\rceil}.$$ Here, $C_{ij} =d_i(\alpha_i,\alpha_j)$ and we say that $j\sim i$ if $C_{i,j}<0$.
Theorem 2.2 of \cite{CMkir} and Theorem 7.1 of \cite{HKOTY} together,  prove the following:
\begin{prop}\label{qsystemg} Assume that $\lie g$ is of type $A,B,C$ or $D$. The $\lie g$--characters of the Kirillov-Reshethikhin modules satisfy the $Q$--system.  More precisely, for $i\in I$ and $m\in\bz_+$, we have a (non--canonical) short exact sequence of $\lie g$--modules, $$0\to K_{i,m}\to KR(m\omega_i)\otimes KR(m\omega_i)\to KR((m+1)\omega_i)\otimes KR((m-1)\omega_i)\to 0,$$ where $$K_{i,m}\cong \bigotimes_{j\sim i}\bigotimes_{k=0}^{-C_{ij}-1} KR(\left({\left\lceil\frac{mC_{ji}-k}{C_{ij}}\right\rceil}\right)\omega_j).$$ \hfill\qedsymbol
\end{prop}

\subsection{} We shall prove the stronger statement:
\begin{thm}\label{qsystem}  Assume that $\lie g$ is of type $A,B,C$ or $D$. Given $i\in I$, and $m\in \bz_+$, there exists a canonical,  non--split short exact sequence of graded $\lie g[t]$--modules:
$$0\to \tau_{d_im} K_{i,d_im}^*\stackrel{\iota}{\to} KR(d_im\omega_i)*KR(d_im\omega_i)\stackrel{\pi}{\to} KR((d_im+1)\omega_i)*KR((d_im-1)\omega_i)\to 0,$$
where $K_{i,d_im}^*$ is obtained from $K_{i,d_im}$ by replacing the tensor product by fusion product and $\tau_{d_im}$ is the grading shift  operator.
\end{thm}

The proof of the theorem occupies the rest of this section. {\em{Moreover $\lie g$ will continue to be  an arbitrary complex simple Lie algebra, unless otherwise stated}}.
\subsection{} Let us first make the tensor product expression for $K_{i,d_im}$ more explicitly. If $d_i=1$ and $j\sim i$, then   $C_{i,j}=-1$, $k=0$, $C_{j,i}=-d_j$ and we have  $$K_{i,d_im}\cong\bigotimes_{j\sim i} KR((d_jm)\omega_j),\ \ d_i=1.$$
If  $d_i=2$ and $j\sim i$ then we have  two possibilities: either $d_j=d_i=2$ in which case $C_{i,j}=-1=C_{j,i}$ and $k=0$ or $d_j=1$ in which case $C_{i,j}=-2$, $C_{j,i}=-1$ and $k=0,1$. Hence we get, $$K_{i, d_im}=\left(\bigotimes_{j\sim i: d_j=2}KR(d_jm\omega_j)\right)\bigotimes\left(\bigotimes_{j\sim i: d_j=1}KR(d_jm\omega_j)^{\otimes 2}\right).$$
If $d_i=3$ then $d_j=1$ and $C_{i,j}= -3$, $C_{j,i}=-1$ and $k\in\{0,1,2\}$. Hence we get $$K_{i,d_im}= KR(m\omega_j)^{\otimes 3}.$$

\subsection{} Set $$\lambda = \begin{cases} m\left(\sum_{j\sim i}d_j\omega_j\right),\ \ \ \ \ d_i=1,\\ \\ m\left(\sum_{j\sim i: d_j=2}d_j\omega_j+ \sum_{j\sim i: d_j=1}2d_j\omega_j\right),\ \ \ \ \ d_i=2,\\ \\ 3m\omega_j,\ \ d_i=3.\\ \end{cases}$$
\vskip 6pt

Using Proposition \ref{krdem}, we see that \begin{equation}\label{kernel} K_{i,d_im}^*\cong_{\lie g[t]} D(m,\lambda)\cong V(\bxi(m,\lambda)),\qquad KR(d_im\omega_i)*KR(d_im\omega_i)\cong V(\bxi(m,2d_im\omega_i)),\end{equation} and hence an equivalent form of Theorem 4 is the following:
\vskip 6pt

\noindent {\bf {Theorem \ref{qsystem}'.}} Assume that $\lie g$ is of type $A$, $B$, $C$ or $D$. Given $i\in I$, there exists a canonical,  non--split short exact sequence of graded $\lie g[t]$--modules:
$$0\to \tau_{d_im} V(\bxi(m,\lambda))\stackrel{\tilde\iota}{\to} V(\bxi(m, 2d_im\omega_i))\stackrel{\tilde\pi}{\to} KR((d_im+1)\omega_i)*KR((d_im-1)\omega_i)\to 0.$$

We now prove Theorem \ref{qsystem}'
  \subsection{} The following Lemma proves the existence of $\tilde \pi$.
   \begin{lem}\label{tau} There exists a surjective map of $\lie g[t]$--modules  $$\tilde\pi: V(\bxi(m,2d_im\omega_i))\to  KR((d_im+1)\omega_i)*KR((d_im-1)\omega_i),$$ such that  $$ \tilde\pi(v_{\bxi(m,2d_im\omega_i)})=w_{i,d_im+1}* w_{i,d_im-1}.$$
   Moreover, $$0\ne (x_i^-\otimes t)^{d_im}v_{\bxi(m,2d_im\omega_i)}\in\ker{\tilde\pi}.$$\end{lem}
   \begin{pf} By Lemma \ref{fusionweyl} we know that   $KR((d_im+1)\omega_i)*KR((d_im-1)\omega_i)$ is a quotient of $W_{\loc}(2d_im\omega_i)$. Hence by Lemma  \ref{first} it is enough to prove that $$(x_i^-\otimes t^2) (w_{i,d_im+1}*w_{i,d_im-1})=0.$$ Let $z_1,z_2$ be the  distinct complex numbers involved in the fusion product. By Lemma \ref{elemfusion}, it is enough to prove that
    $$(x_i^-\otimes (t-z_1)(t-z_2))(w_{i,d_im+1}*w_{i,d_im-1})=0.$$ This follows from,
\begin{gather*} (x_i^-\otimes (t-z_1)(t-z_2))(w_{i,d_im+1}\otimes w_{i,d_im-1})=\\ (x_i^-\otimes t(t-z_2+z_1))w_{i,d_im+1}\otimes  w_{i,d_im-1}+ w_{i,d_im+1}\otimes  (x_i^-\otimes t(t-z_1+z_2))w_{i,d_im-1} = 0,\end{gather*} where the last equality is a consequence of equation \eqref{bigger}.  This proves  the existence of $\tilde\pi$.
\vskip6pt

 We now prove that  $$(x_i^-\otimes t)^{d_im}(w_{i,d_im+1}* w_{i,d_im-1})=0.$$ Again, using Lemma \ref{elemfusion}, we have,  $$(x_i^-\otimes t)^{d_im}(w_{i,d_im+1}* w_{i,d_im-1})=(x_i^-\otimes (t-z_1))^{d_im}(w_{i,d_im+1}* w_{i,d_im-1}), $$ and this time we get,  \begin{gather*} (x_i^-\otimes (t-z_1))^{d_im}(w_{i,d_im+1}\otimes  w_{i,d_im-1})= \\ \left((x_i^-\otimes t)^{d_im}w_{i,d_im+1}\right)\otimes  w_{i,d_im-1} + w_{i,d_im+1}\otimes  (x_i^-\otimes (t-z_1+z_2))^{d_im}w_{i,d_im-1}.\end{gather*}The first term on the right hand side is zero by equation \eqref{bigger} while the second term is zero because $$(d_im-1)\omega_i-(d_im\alpha_i)\notin\wt(KR((d_im-1)\omega_i)).$$ This proves,  that $$(x_i^-\otimes t)^{d_im}v_{\bxi(m,2d_im\omega_i)}\in\ker{\tilde\pi}.$$

  \vskip6pt

\noindent  To complete the proof of the Lemma, we must show
that $$(x_i^-\otimes t)^{d_im}v_{\bxi(m,2d_im\omega_i)}\ne 0, $$ and using the isomorphism $$V(\bxi(m,2md_i\omega_i))\cong KR(d_im\omega_i)*KR(d_im\omega_i),$$  it is enough to show that $\bov= (x_i^-\otimes t)^{d_im}(w_{i,d_im}*w_{i,d_im})$ is a non--zero element  of the fusion product $KR(d_im\omega_i)*KR(d_im\omega_i)$.  Consider  the tensor product $$\bk= KR(d_im\omega_i)^{z_1}\otimes KR(d_im\omega_i)^{z_2}.$$ Then, we have
\begin{gather*}(x_i^-\otimes (t-z_1))^{d_i,m}(w_{i,d_im}\otimes w_{i,d_im})= w_{i,d_im}\otimes  (x_i^-\otimes (t-z_1+z_2))^{d_im}w_{i,d_im} \\  =  (z_2-z_1)^{d_im}w_{i,d_im}\otimes (x_i^-)^{d_im}w_{i,d_im}.\end{gather*}  To prove that $\bov\ne 0$ it suffices from the preceding equalities to prove that $$\bow= w_{i,d_im}\otimes (x_i^-)^{d_im}w_{i,d_im}\notin\cal F^r,  \ \ r<d_im,$$ where  $\cal F^r$ is the $r$--th filtered piece of $\bk$. Since $\bow\in\bk_{d_im(\omega_i-\alpha_i)}$, this amounts to proving that $$\bow\notin{\rm{span}}\{(x_i^-\otimes 1)^s(x_i^-\otimes (t-z_1))^{d_im-s}(w_{i,m}\otimes w_{i,m}): s> 0\},$$ or equivalently,
$$\bow\notin{\rm{span}}\left\{(x_i^-\otimes 1)^s\left(w_{i,m}\otimes (x_i^-\otimes 1)^{d_im-s}w_{i,m}\right): s> 0\right\}.$$ But this is now a simple calculation using the comultiplication for $\lie g$.

   \end{pf}

 \subsection{}  The next result establishes the existence of $\tilde \iota$. \begin{lem}\label{iota} There exists an injective  non--zero map of graded $\lie g[t]$--modules, $$\iota:\tau_{d_im}K_{i,d_im}^*\to \ker\pi.$$
  \end{lem}
  \begin{pf}  Since $K_{i,d_im}^*\cong V(\bxi(m,\lambda))$, the Lemma will follow if we prove that the non--zero element $(x_i^-\otimes t)^{d_im}v_{\bxi(m,2d_im\omega_i)}\in\ker \tilde\pi$ satisfies the relations of $\tau_{d_im}V(\bxi(m,\lambda))$. It is trivial to see that the element satisfies the relations of the grade shifted  local Weyl module $\tau_{d_im}W_{\loc}(\lambda)$. By Lemma \ref{first}  we only need to check that the following relations hold for all $j\sim i$:
 \begin{gather*}(x_j^-\otimes t)(x_i^-\otimes t)^{d_im}v_{\bxi(m,2d_im\omega_i)}= 0,\ \ \ \ d_i=1\ \ {\rm{or}}\ \ d_i=d_j=2,\\
  (x_j^-\otimes t^2)(x_i^-\otimes t)^{d_im}v_{\bxi(m,2d_im\omega_i)}=0,\ \ d_i=2,\ \  d_j=1.\end{gather*}  If $i=j$, the result is immediate from \eqref{demrelo}. If $i\ne j$, then  $\alpha_i,\alpha_j$ span a root system of type $A_2$ or $C_2$, we see, by using the commutation relations in these Lie algebras, that the equalities follow if we prove that \begin{equation}\label{firstkind} (x^-_{\alpha_i+\alpha_j}\otimes t^2)v_{\bxi(m,2d_im\omega_i)}= 0,\end{equation} and also if $d_i=2$, $d_j=1$, then \begin{equation}\label{secondkind}(x^-_{2\alpha_i+\alpha_j}\otimes t^4)v_{\bxi(m,2d_im\omega_i)}=0. \end{equation} These  are  proved by a straightforward case by case analysis using the relations \eqref{demrelo}. Thus we have,
  \begin{gather*} d_i=d_j\implies d_{\alpha_i+\alpha_j}=d_i,\ \  \omega_i(h_{\alpha_i+\alpha_j})=1,\\
d_i=1,\ d_j=2\implies d_{\alpha_i+\alpha_j}=2,\ \ \omega_i(h_{\alpha_i+\alpha_j})=2,\\
d_i=2,\  \ d_j=1\implies d_{\alpha_i+\alpha_j}=2,\ \ \omega_i(h_{\alpha_i+\alpha_j})= 1,\\
d_i=2, d_j=1\implies d_{2\alpha_i+\alpha_j}=1,\ \  \omega(h_{2\alpha_i+\alpha_j})=1.
\end{gather*}
In the first three  cases, this means that $2d_i m\omega_i(h_{\alpha_i+\alpha_j})=2d_{\alpha_i+\alpha_j} m$ and in the last case, we have $2d_im\omega_i(h_{2\alpha_i+\alpha_j})= 4d_{2\alpha_i+\alpha_j}m$. Equations \eqref{firstkind} and \eqref{secondkind} are now immediate from \eqref{demrelo}. The case when $d_i=3$ is similar and we omit the details.
Since $\tilde\iota$ is a non--zero map between Demazure modules, it follows  from Theorem \ref{demprop}(i) that $\tilde\iota$ is injective and the proof of the Lemma is complete.

  \end{pf}
  \subsection{}
   Note that Lemma \ref{tau} and Lemma \ref{iota} prove that
  $$\dim\left(KR(d_im\omega_i)*KR(d_im\omega_i)\right)\ge\dim K_{d_im}^*+\dim\left( KR((d_im+1)\omega_i)*KR((d_im-1)\omega_i)\right).$$ Hence Theorem \ref{qsystem} obviously  follows if we establish the reverse inequality.
  Since the dimension is unchanged when we replace tensor products by fusion products, we must prove that \begin{equation}\label{dimequal}\dim\left(KR(d_im\omega_i)\otimes KR(d_im\omega_i)\right)= \dim K_{d_im}+\dim\left( KR((d_im+1)\omega_i)\otimes KR((d_im-1)\omega_i)\right).\end{equation} But this is immediate from Proposition \ref{qsystemg} if $\lie g$ is of classical type. Finally, since
   $V(\bxi(m,2d_im\omega_i))$ is an indecomposable module, it follows that the sequence is non--split and the theorem is proved.

\subsection{}\label{general} Theorem \ref{qsystem} can also be deduced for arbitrary simple Lie algebras from Lemma \ref{tau} and Lemma \ref{iota} in a different way and we
 explain this very briefly.   The Kirillov--Reshetikhin modules for the current algebra were proved in \cite{Ckir}, \cite{kedem} to be the specializations to $q=1$ of Kirillov--Reshetikhin modules for the quantized enveloping algebra of the loop algebra. An obvious consequence of results proved in \cite{H1}, \cite{NakT1},\cite{NakT2}  is that the quantum Kirillov--Reshetikhin modules satisfy the required conditions on the dimensions in \eqref{dimequal}. Since the dimension is unchanged on specializing we have the necessary equality of dimensions for the current algebra. However all these results are very difficult, use many different ideas and deep results in the representation theory of quantum and classical affine Lie algebras.
 Our approach keeps us in the purely classical realm of  Lie algebras and shows that the only missing piece of information for proving Theorem \ref{qsystem} for exceptional Lie algebras is \eqref{dimequal}. Since the Kirillov--Reshetikhin modules are actually Demazure modules, it might be possible to use the literature on Demazure modules to prove \eqref{dimequal} directly, although this too seems hard at the moment.

\section{The $\lie{sl}_2$ case}
In this section, we focus on the Lie algebra $\lie{sl}_2$ and prove that the modules $V(\bxi)$ are fusion products of evaluation modules $\ev_0V(r)$, $r\in\bz_+$ and vice--versa. In particular, this gives generators and relations for such fusion products and recovers in a very different and self contained way the results of \cite{FF}.
We also give an explicit monomial basis for the modules $V(\bxi)$.

\subsection{} In the case of $\lie{sl}_2$, the set $I$ is the singleton set $\{1\}$
and so we just denote the elements $x_1^\pm$ by $x^\pm$ etc.
We  identify $P$ with the integers and given $r\in\bz_+$, a $r$--compatible partition $\xi$ is just a partition of $r$, i.e.,  $\xi=(\xi_1\ge\xi_2\ge \cdots ) $ and  $|\xi|=r$.

\subsection{}\label{ximp} Given a partition  $\xi=(\xi_1\ge\cdots\ge \xi_\ell>0)$ with $\ell$ parts we define  partitions $\xi^\pm$ as follows. If $\ell=1$, then  $\xi^+=\xi$ and $\xi^-$ is the empty partition. If $\ell>1$, then  $\xi^-=(\xi^-_1\ge\cdots \ge\xi^-_{\ell-2}\ge\xi^-_{\ell-1}\ge 0)$ is given by  $$\xi^-_r=\begin{cases}\xi_r,\ \ r<\ell-1,\\ \xi_{\ell-1}-\xi_\ell,\ \ r=\ell-1,\\ 0,\  \ r\ge \ell.\end{cases}$$ In particular $|\xi^-|=|\xi|-2\xi_\ell$.
We take  $\xi^+=(\xi_1^+\ge\cdots\ge \xi^+_{\ell-1}\ge\xi_\ell^+\ge 0)$ is the unique  partition associated to the $n$--tuple $(\xi_1,\cdots, \xi_{\ell-2},\xi_{\ell-1}+1,\xi_{\ell}-1)$. It can be described explicitly as follows:  let  $0\le \ell(\xi) \le  \ell-2$ be minimal such that $\xi_{\ell(\xi)+1}=\xi_{\ell-1}$.

Then, $$\xi^+_j=\begin{cases}\xi_j,\ \ \ \  1\le  j\le \ell(\xi)\ \ {\rm{or}}\ \ j>\ell, \\
\xi_{\ell-1}+1,\ \ \ \ j=\ell(\xi)+1,\\
\xi_{\ell-1}, \ \ \ \ \ell(\xi)+2\le j\le \ell-1,\\
\xi_\ell-1,\ \ \ \  j=\ell.\end{cases}$$

and hence,
 \begin{equation}\label{relxi}\sum_{j\ge k+1}\xi_j^+=\begin{cases} \sum_{j\ge k+1}\xi_j,\  \ \ \ k\le\ell(\xi)\ \ {\rm{or}}\ \ k\ge \ell,\\ \\ -1+ \sum_{j\ge k+1}\xi_j,\ \ \  \ \ell(\xi)<k\le \ell-1.\\
 \end{cases}\end{equation}

 \subsection{} Define a subset  $\mathbb I(\xi)$  of $\bz_+^\ell$ by:  $\boi=(i_1,\cdots,i_\ell)\in\mathbb I(\xi)$, iff for all $2\le k\le \ell+1$ and $1\le j\le k-1$, we have,
$$(ji_{k-1}+ (j+1)i_k)+2\sum_{p=k+1}^\ell i_p\le \sum_{p=k-j}^\ell \xi_p,$$ where we understand that $i_{\ell+1}=0$. We shall prove,
\begin{thm}\label{sl2} Assume that $\lie g$ is isomorphic to $\lie{sl}_2$ and that $\xi=(\xi_1\ge\cdots\ge \xi_\ell>0)$.
\begin{enumerit}
\item[$(i)$] For  $\ell>1$, there exists a short exact sequence of $\lie g[t]$--modules,$$0\to\tau_{(\ell-1)\xi_\ell}V(\xi^-)\stackrel{\varphi^-}{\longrightarrow} V(\xi){\stackrel{\varphi^+}{\longrightarrow}} V(\xi^+)\to 0.$$
\item[$(ii)$] For all $\ell\ge 1$, we have an isomorphism of $\lie g[t]$--modules, $$V(\xi)\cong V^{z_1}(\xi_1)*\cdots *V^{z_\ell}(\xi_\ell),$$ for any set $z_1,\cdots, z_\ell$ of distinct scalars.
\item[(iii)] The elements $$\{{(x^-\otimes 1})^{i_1}\cdots (x^-\otimes t^{\ell-1})^{i_\ell}v_\xi: (i_1,\cdots, i_\ell)\in\mathbb I(\xi)\},$$ form a basis of $V(\xi)$.\end{enumerit}
\end{thm}

\begin{rem}  The existence of this exact sequence was proved in \cite{FF}, \cite{FF2} by very different methods and do not use the definition of the fusion product in an essential way. We are able to give very direct proofs of these statements because of the three presentations of the modules $V(\xi)$, established earlier in the paper. \end{rem}
The proof of the Theorem occupies the rest of this section. The first two parts of the theorem are proved simultaneously by an induction on $\ell$, namely we shall prove that if $(ii)$ holds for a partition with $\ell$ parts, then $(i)$ and $(ii)$  hold for a partition with $\ell+1$ parts. Proposition \ref{ev0} shows that $(ii)$ holds when $\ell=1$ and so induction begins.
\noindent {\em{ In the rest of the section, given $s,r\in\bz_+$, we set, $$X(r,s)=(x^+\otimes t)^{(s)}(x^-\otimes 1)^{(s+r)}.$$}}
\subsection{}\label{pi1} The following result establishes the existence of $\varphi^+$.
\begin{prop}  There exists a surjective morphism of $\lie g[t]$--modules $\varphi^+:V(\xi)\to V(\xi^+)$ such that $$\varphi(v_\xi)=  v_{\xi^+},\ \qquad  \ker\varphi^+=\bu(\lie g[t])(x^-\otimes t^{\ell-1})^{\xi_\ell}v_\xi.$$

\end{prop}
\begin{pf} To prove that $\varphi^+$ exists, it suffices to show  that $$X(r,s)v_{\xi^+}=0,\ \  {\rm{for\ all}}\ \ s, r,k\in\bn, \ {\rm{with}}\  s+r\ge 1+ rk+\sum_{j\ge k+1}\xi_j.$$ But this is immediate by noting that $\sum_{j\ge k+1}\xi_j\ge\sum_{j\ge k+1}\xi^+_j$.
\vskip 6pt

\noindent Using \eqref{demreltype} we see  that $$(x^-\otimes t^{\ell-1})^{\xi_\ell}v_{\xi^+}=0,\qquad {\rm{i.e.}},\ (x^-\otimes t^{\ell-1})^{\xi_\ell}v_\xi\in\ker\varphi^+.$$ To complete the proof we must show that  $(x^-\otimes t^{\ell-1})^{\xi_\ell}v_\xi$ generates $\ker\varphi$. Suppose that $Xv_{\xi^+}=0$ for some  $X\in\bu(\lie g[t])$. Then we may write $X=Y+Z$ where $Y$ is in the left ideal of $\bu(\lie g[t])$ generated by the set $$\{ x^+\otimes t^p, \ \  (h\otimes t^p)-\delta_{p,0}|\xi|, \ \ (x^-\otimes 1)^{|\xi|+1}: p\in\bz_+\},$$ and $Z$ is in the left ideal of $\bu(\lie g[t])$ generated by
  $$\{X(r,s): s,r,k\in\bn, \ \ s+r\ge 1+rk+\sum_{j\ge k+1}\xi^+_j\ \}.$$ Since $Yv_\xi=0$ as well, we need to consider only the case when $X=X(r,s)$ for $s,r\in\bn$ with  $s+r\ge 1+rk+\sum_{j\ge k+1}\xi_j^+$ for some $k\in\bn$. Moreover, we see by using equation \eqref{relxi}, that  the only time $X(r,s)v_\xi=0$  is not obviously a defining relation of $v_\xi$ is when $s+r= rk+\sum_{j\ge k+1}\xi_j$,\ \ {\rm{for\ some}}\ \ $\ell(\xi)<k\le \ell-1$. Since $\xi_j=\xi_{\ell-1}$ if $\ell(\xi)+2\le j\le \ell- 1$ we have $$\ell(\xi)<k\le\ell-1\implies \ \ s+r=rk+(\ell-k-1)\xi_{\ell-1}+\xi_\ell.$$
 If $\xi_{\ell-1}<r$, then we
have   $$s+r=(k-1)r+ (r-\xi_{\ell-1})+\xi_{\ell-1}(\ell-k)+\xi_\ell\ge 1+  (k-1)r+\sum_{j\ge k}\xi_j,$$
(since $\xi_k=\xi_{\ell-1}$ under the hypothesis that $\ell-1\ge k>\ell(\xi)$).  This means that  $X(r,s) v_\xi=0$ and we are done in this case.

Assume that $\xi_{\ell-1}\ge r$ and let  $${(b_p)_{p\ge 0}\in\bs(r,s)},\qquad s+r=rk+(\ell-k-1)\xi_{\ell-1}+\xi_\ell$$ for some $\ell(\xi)<k\le \ell-1$. Recalling from Corollary \ref{third} that
 $(x^-\otimes t^\ell)v_\xi=0$, we see that $${(x^-\otimes 1)}^{(b_0)}{(x^-\otimes t)}^{(b_1)}\cdots {(x^-\otimes t^{s})}^{(b_{s})}v_\xi\ne 0 \implies s\le (\ell-1)b_{\ell-1}+(r-b_{\ell-1})(\ell-2).$$ This gives,
 $$r(k-1)+(\ell-k-1)\xi_{\ell-1}+\xi_\ell\le b_{\ell-1}+r(\ell-2),$$ i.e., $$0\le (\ell-k-1)(\xi_{\ell-1}-r)\le b_{\ell-1}-\xi_\ell.$$  Hence $b_{\ell-1}\ge\xi_\ell$ and so, $$X(r,s)v_\xi= \bold{x}^-(r,s)v_\xi\in \bu(\lie n^-[t])(x^-\otimes t^{\ell-1})^{\xi_\ell}v_\xi, $$
which  completes the proof.
\end{pf}

\subsection{} We need some additional results to   prove the existence of $\varphi^-$. The following Lemma is checked by a straightforward induction on $p$.
\begin{lem}\label{simple} Given $a,b,p\in\bz_+$, we have,
\begin{gather*}   [(x^+\otimes t^a)\ ,\  (x^-\otimes t^b)^{(p)}]= (x^-\otimes t^b)^{(p-1)}(h\otimes t^{a+b})-(x^-\otimes t^{a+2b})(x^-\otimes t^b)^{(p-2)},\\
[(h\otimes t^a) \ ,\ (x^-\otimes t^b)^{(p)}] =-2(x^-\otimes t^{a+b})(x^-\otimes t^b)^{(p-1)},\\
\label{one}[(x^-\otimes t^{a})\ ,\ (x^+\otimes t^{b})^{(p)}]= -(x^+\otimes t^{b})^{(p-1)}(h\otimes t^{a+b})-(x^+\otimes t^b)^{(p-2)}(x^+\otimes t^{a+2b}),\end{gather*}
where we understand as usual that $(x\otimes t^b)^p=0$ if $p<0$.\hfill\qedsymbol

\end{lem}
\subsection{} The next result we need is the following.
\begin{lem}\label{critical} Suppose that $V$ is a $\lie g[t]$--module and   $v\in V$ satisfies,
\begin{gather*}\label{lweyl}\lie n[t]v =0, \quad  \ \lie h[t]_+v=0, \ \quad\ \ X(r,s)v=0,  \ \ \  s,r\in\bn \ \ s+r\ge N,\end{gather*} for some $N\in\bn$.
 Then for all $j\in\bn$, we have  $$X(r,s)(x^-\otimes t^{j})v=0  \ {\rm{for\ all}}\ \  s,r\in\bz_+\ {\rm{with}}\ \ s+r\ge N-2.$$
\end{lem}
\begin{pf} For $j\in \bn$, set $v_j=(x^-\otimes t^{j})v$ and note that $\lie n^+[t]v_j=0$. Write $$(s+1)X(r+1,s+1)=(x^+\otimes t)^{(s)}(x^+\otimes t)(x^-\otimes 1)^{(r+s+2)}.$$
Using the first equation in Lemma \ref{simple} and the relations satisfied by $v$, we find that $$(s+1)X(r+1,s+1)v= -X(r,s)v_1.$$
 Since the left hand side is zero if $s+r+2\ge N$, it follows that \begin{equation}\label{one}X(r,s)v_1=0,\ s+r\ge  N-2.\end{equation}
 Using the second and third equations in Lemma \ref{simple} we see that, for $i\ge 0$, we have \begin{gather*} (x^-\otimes t^{i})X(r,s)v=
X(r,s)(x^-\otimes t^i)v+2X(r,s-1)v_{i+1}+ X(r,s-2)v_{i+2}.\end{gather*}
If we take $s,r$ such that $s+r\ge N$, then the left hand side of the equation is zero, i.e.,$$ X(r,s)(x^-\otimes t^i)v+2X(r,s-1)v_{i+1}+ X(r,s-2)v_{i+2}=0,\ \ s+r\ge N.$$ If $i=0$, then the first term on the left hand side is a multiple of $X(r+1,s)v$ and hence is zero.  Using \eqref{one} we see that since $s+r-1\ge N-2$ the second term is zero which proves finally that $$X(r,s-2)v_{2}=0, \ s+r\ge N,\ \ {\rm{i.e.}},\ \  X(r,s)v_2=0,\ \ s+r\ge N-2.$$ An obvious induction on $i$ now gives the lemma.

\end{pf}
\begin{cor} Let $v\in V$ be as in the lemma and assume also that $(x^-\otimes t^p)v=0$. Then $$X(r,s)(x^-\otimes t^{p-1})^kv=0,\ \ s,r\in\bz_+,\ \  s+r\ge N-2k.$$
\end{cor}
\begin{pf} The proof is immediate once we observe that the element $(x^-\otimes t^{p-1})^kv$ satisfies the same relations as $v$ with $N$ replaced by $N-2k$.
\end{pf}

\subsection{}\label{injectivesl2} We now establish the existence of $\varphi^-$.
\begin{prop}

There exists a surjective morphism of $\lie g[t]$--modules satisfying
$$\varphi^-:V(\xi^-)\to \ker\varphi^+,\qquad\varphi^-(v_{\xi^-})= (x^-\otimes t^{\ell-1})^{\xi_\ell}v_\xi.$$
\end{prop}
\begin{pf} It is  trivial to check that there exists a morphism $W_{\loc}(|\xi^-|)\to V(\xi)$ which maps $w_{|\xi^-|}\to (x^-\otimes t^{\ell-1})^{\xi_\ell}v_\xi.$
 To prove that the map factors through to  $V(\xi^-)$  we must show that
 $$X(r,s)(x^-\otimes t^{\ell-1})^{\xi_\ell}v_\xi=0,\ \ s, r,k\in\bn, \ {\rm{with}}\  s+r\ge 1+ rk+\sum_{j\ge k+1}\xi_j^-.$$

  If $k\le \ell-2$, then we have $$1+ rk+\sum_{j\ge k+1}\xi_j^-=1+ rk+\left(\sum_{j\ge k+1}\xi_j\right)-2\xi_\ell.$$ Since $X(r,s)v_\xi=0$ if $s+r\ge 1+\sum_{j\ge k+1}\xi_j$ and $$(x^-\otimes t^\ell)v_\xi=0\implies (x^-\otimes t^\ell)(x^-\otimes t^{\ell-1})^{\xi_\ell}v_\xi=0,$$ the result is immediate from Corollary  \ref{critical}.

  If $k\ge\ell-1$, then $s\ge 1+r(\ell-2)$ and hence  $$(b_p)_{p\ge 0}\in \bs(r,s) \implies\  b_p>0,\ \ {\rm{for\ some}}\ \  p\ge \ell-1.$$ This means that $$\bx^-(r,s) =\sum_{p\ge\ell-1} X_p(x^-\otimes  t^p)$$ for some $X_p\in\bu(\lie n^-[t])$ which gives, $$\bx^-(r,s)(x^-\otimes t^{\ell-1})^{\xi_\ell}v_\xi= (x^-\otimes t^{\ell-1})^{\xi_\ell}\left(\sum_{p\ge\ell-1} X_p(x^-\otimes  t^p)\right)v_\xi.$$
Now $(x^-\otimes t^{\ell-1})^{\xi_\ell+1}v_{\xi}=0$ by  proposition \ref{third}, so equivalently we get $X(r,s)(x^-\otimes t^{\ell-1})^{\xi_\ell}v_{\xi}=0$ as required.

\end{pf}

\subsection{} The existence of the surjective map $\varphi^+$ and the map $\varphi^-$ implies that we have \begin{equation}\label{dimsmaller}\dim V(\xi)\le\dim V(\xi^+)+\dim V(\xi^-).\end{equation} To prove the reverse inequality we need the following.
\begin{prop}\label{dimbigger} The assignment $v_\xi\to v_{\xi_1}*\cdots *v_{\xi_\ell}$ defines a surjective homomorphism $$V(\xi)\to V(\xi_1)*\cdots* V(\xi_\ell)\to 0.$$ In particular, $$\dim V(\xi)\ge \dim V(\xi_1)\cdots\dim V(\xi_\ell)=\prod_{s=1}^\ell(\xi_s+1).$$ \end{prop}
Assuming the proposition, the proof of parts (i) and (ii) of the Theorem \ref{sl2} is completed as follows. Assume that $(i)$ and $(ii)$ hold for all $1\le m\le\ell-1$. We must prove that $(i)$ and $(ii)$ hold for $\ell$.     If $\xi$ has $\ell$ parts and $\xi_\ell=1$, then $\xi^\pm$ have at most $(\ell-1)$ parts and hence by induction we see that $$\dim V(\xi)\ge\prod_{s=1}^\ell(\xi_s+1)= \prod_{s=1}^{\ell-1}(\xi^+_s+1)+\prod_{s=1}^{\ell-1}(\xi^-_s+1)=\dim V(\xi^+)+\dim V(\xi^-).$$ Together with \eqref{dimsmaller} we get $$\dim V(\xi)=\prod_{s=1}^\ell(\xi_s+1)=\dim V(\xi^+)+ \dim V(\xi^-),$$ and parts (i) and (ii) are proved in this case.  An obvious further  induction on $\xi_\ell$ now completes the proof in the case when $\xi$ has $\ell$ parts and proves the inductive step.

\subsection{} We need some additional comments to prove Proposition \ref{dimbigger}. The following Lemma  is trivially checked.
\begin{lem}\label{localize} Given any $f\in\bc[t]$, we have an injective map $\psi_f: \lie{sl}_2\otimes\bc[t]\to\lie{sl}_2\otimes\bc[t,1/f]$ of Lie algebras  \begin{gather*}\psi_f(x^+\otimes t^p)= x^+\otimes t^p/f,\ \ \ \ \psi_f(x^-\otimes t^p)= x^- \otimes t^pf,\ \ \ \psi(h\otimes t^p)= h\otimes t^p.\end{gather*} In particular, if $a\in\bc$ is such that $f(a)\ne 0$, and $V$ is any representation of $\lie g$, then $\ev_aV$ is a representation of $\lie{sl}_2\otimes\bc[t,1/f]$, given by $$(x\otimes g)v=g(a)xv,\ \ g\in\bc[t,1/f], \ \ x\in\lie g.$$\hfill\qedsymbol
\end{lem}
\noindent We  note some consequences of the Lemma. For $s,r,k\in\bn$, we have \begin{equation}\label{shift}\psi_{t^k}(\bx^-(r,s))=\ \ _k\bx^-(r,s+rk).\end{equation}

\noindent For $s,r\in\bn$ and $f\in\bc[t]$, we have, by using Lemma \ref{gar}, \begin{equation}\label{crit}\psi_f(X(r,s))= (-1)^s\psi_f(\bx^-(r,s)) +\sum_{p\ge 1}\psi_f(Y_p)H_p+\psi_f(X),\end{equation} for some $Y_p\in\bu(\lie n^-[t])[s-p]$,  $H_p\in\bu(\lie h[t]_+)[p]$ and $X\in\bu(\lie g[t])(\lie n^+[t])$.

Suppose that $f\in\bc[t]$ and  $a_1,\cdots,a_m\in\bc$,  $m\in\bn$  are such that $f(a_s)\ne 0$ for $1\le s\le m$.  If $V_1,\cdots ,V_m$ are arbitrary modules for $\lie{sl}_2$, then $$V=\ev_{a_1}V_1\otimes\cdots \otimes \ev_{a_m}V_m$$ is a module for $\lie{sl}_2[t,1/f]$.
\vskip 6pt

\subsection{}{\em {Proof of Proposition \ref{dimbigger}.}} Using Lemma \ref{fusionweyl} and  Proposition \ref{third}, it is enough to prove that $$_{k}\bx^-(r,s)(v_{\xi_1}*\cdots*v_{\xi_\ell})=0,\ \ s,r,k\in\bn, s+r\ge 1+kr+\sum_{j\ge k+1}\xi_j.$$ Using \eqref{shift} this is equivalent to proving that $$\psi_{t^k}(\bx^-(r,s-rk))(v_{\xi_1}*\cdots*v_{\xi_\ell})=0,\ \ s,r,k\in\bn, s+r\ge 1+kr+\sum_{j\ge k+1}\xi_j.$$ Let $z_1,\cdots,z_\ell$ be the distinct parameters used to define the fusion product. For $1\le k\le \ell$, set  $f_k=(t-z_1)\cdots (t-z_{k})$. Using Lemma \ref{elemfusion} we see that $$\psi_{t^k}(\bx^-(r,s-rk))(v_{\xi_1}*\cdots*v_{\xi_\ell})=\psi_{f_k}(\bx^-(r,s-rk))(v_{\xi_1}*\cdots*v_{\xi_\ell}). $$
Since $$(x\otimes f_k)(V(\xi_1)\otimes \cdots\otimes V(\xi_k))=0,\ \ x\in\lie g,$$ and $\psi_{f_k}(x^-\otimes t^p)=(x^-\otimes t^pf_k)$, we see that,
$$\psi_{f_k}(\bx^-(r,s-rk))(v_{\xi_1}\otimes\cdots\otimes v_{\xi_\ell})=(v_{\xi_1}\otimes\cdots\otimes v_{\xi_k})\otimes\psi_{f_k}(\bx^-(r,s-rk))(v_{\xi_{k+1}}\otimes\cdots\otimes v_{\xi_\ell}).$$   Using \eqref{crit} and the fact that $$(\lie n^+[t,1/f_k])(v_{\xi_{k+1}}\otimes\cdots\otimes v_{\xi_\ell})=0,$$ we get,\begin{gather*}((-1)^{s}\psi_{f_k}(\bx^-(r,s-rk)) +\sum_{p\ge 1}\psi_{f_k}(Y_p)H_p-\psi_{f_k}(X(r,s-rk)))
(v_{\xi_{k+1}}\otimes\cdots\otimes v_{\xi_\ell})=0,\end{gather*}
 where $Y_p\in\bu(\lie n^-[t])[s-rk-p]$, and $H_p\in\bu(\lie h[t]_+)[p]$.
 Now,
\begin{gather*} \psi_{f_k}((x^+\otimes t)^{s-rk}(x^-\otimes 1)^{s-r(k-1)})(v_{\xi_{k+1}}\otimes\cdots\otimes v_{\xi_\ell})=0,\end{gather*} since $s-r(k-1)\ge 1+\sum_{j\ge k+1}\xi_j$ and so cannot be a weight of $V(\xi_{k+1})\otimes\cdots\otimes V(\xi_\ell)$, i.e.$$\psi_{f_k}(X(r,s-rk))(v_{\xi_{k+1}}\otimes\cdots\otimes v_{\xi_\ell})= 0,\ \ s+r\ge 1+rk+\sum_{j\ge k+1}\xi_j.$$
Further $$H_p(v_{\xi_{k+1}}\otimes\cdots\otimes v_{\xi_\ell})=A_p(v_{\xi_{k+1}}\otimes\cdots\otimes v_{\xi_\ell}),$$ for some $A_p\in\bc$ and we have, \begin{gather*}(-1)^s\psi_{f_k}(\bx^-(r,s-rk))(v_{\xi_{k+1}}\otimes\cdots\otimes v_{\xi_\ell})=\sum_{p\ge 1}A_p\psi_{f_k}(Y_p)(v_{\xi_{k+1}}\otimes\cdots\otimes v_{\xi_\ell}).\end{gather*}
 Putting all this together we find finally that
$$\psi_{f_k}(\bx^-(r,s-rk))(v_{\xi_1}\otimes\cdots\otimes v_{\xi_\ell})=(-1)^s\sum_{p\ge 1}A_p\psi_{f_k}(Y_p)(v_{\xi_1}\otimes\cdots\otimes v_{\xi_k}\otimes v_{\xi_{k+1}}\otimes\cdots\otimes v_{\xi_\ell}).$$ Since $\psi_{f_k}(Y_p)\in\bu(\lie n^-[t])[s-p]$ and $p\ge 1$,  we get that $$\psi_{f_k}(\bx^-(r,s-rk))(v_{\xi_1}*\cdots* v_{\xi_\ell}))=0,$$ as needed.

It remains to consider the case when $k\ge \ell+1$. But here we have $$\psi_{f_\ell}(\bx^-(r,s))(v_{\xi_1}\otimes\cdots\otimes v_{\xi_\ell})=0,$$ and the proof is complete.

\subsection{}\label{basissl2} We shall deduce Theorem \ref{sl2}(iii) from the following proposition.
\begin{prop} We have,  $$\mathbb I(\xi)=\mathbb I(\xi^-;\xi_\ell)\ \bigsqcup\ \mathbb I(\xi^+)$$ where $$ \mathbb I(\xi^-;\xi_\ell)=\{(i_1,\cdots, i_{\ell-1},\xi_\ell):(i_1,\cdots, i_{\ell-1})\in \mathbb I(\xi^-)\}.$$
\end{prop}
\begin{pf} Using the explicit formulae for $\xi^+$  given in Section \ref{ximp}, it is trivial to see that  $$\mathbb I(\xi^+)\subset\mathbb I(\xi).$$ Note also that $$(i_1,\cdots, i_\ell)\in\mathbb I(\xi^+)\implies i_\ell<\xi_\ell.$$
To see this, recall from the definition of $I(\xi^+)$ that the inequality $$i_\ell\le\xi^+_\ell=\xi_\ell-1,$$ must be satisfied.
\vskip 6pt

Suppose that $\boi=(i_1,\cdots, i_{\ell-1},\xi_\ell)\in \mathbb I(\xi^-;\xi_\ell)$, in which case, for all $2\le k\le \ell$ and $1\le j\le k-1$, we have, $$ji_{k-1}+ (j+1)i_k+2(i_{k+1}+\cdots+ i_{\ell-1})\le  \xi_{\ell-1}-\xi_\ell+ \sum_{p=k-j}^{\ell-2} \xi_p$$ which is clearly equivalent to $$ji_{k-1}+ (j+1)i_k+2(i_{k+1}+\cdots +i_{\ell-1}+\xi_\ell)\le   \sum_{p=k-j}^{\ell} \xi_p.$$  If $k=\ell+1$ and $1\le j\le \ell$, then we  have $$j\xi_\ell\le \sum_{p=\ell+1-j}^\ell\xi_p, $$ and hence we have proved that $$\mathbb I(\xi^-;\xi_\ell)\subset \mathbb I(\xi),\qquad \mathbb I(\xi^+)\cap\mathbb I(\xi^-;\xi_\ell)=\emptyset.$$
\vskip 6pt
\noindent Notice also that the preceding argument also proves that\begin{equation*}\boi\in\mathbb I(\xi),\ \  i_\ell=\xi_\ell\implies (i_1,\cdots, i_{\ell-1})\in \mathbb I(\xi^-),\end{equation*} i.e. $$I(\xi^-;\xi_\ell)=\{(i_1,\cdots,i_\ell)\in\mathbb I(\xi): \ \ i_\ell=\xi_\ell\}.$$ Thus, to complete the proof of the proposition
we must prove,\begin{equation}\label{final}(i_1,\cdots, i_\ell)\in\mathbb I(\xi),\ \ i_\ell< \xi_\ell\implies \boi\in\mathbb I(\xi^+).\end{equation} \iffalse Fix $\ell$ minimal such that $\xi_\ell=\xi_n$, i.e. $\xi_{r}=\xi_n$ if $r\ge \ell$ and $\xi_{\ell-1}>\xi_n$. Then,
$$\xi^+=\xi_1\ge\xi_2\ge\cdots\ge \xi_{\ell-1}\ge\xi_{n}+1>\xi_n\ge\cdots\ge\xi_n>\xi_n-1.$$ \fi
If $2\le k\le \ell+1$ and $1\le j\le k-1$ are such that $k-j\le\ell(\xi)+1$, then the inequalities in  $\mathbb I(\xi)$ and $\mathbb I(\xi^+)$ are the same and there is nothing to prove. If  $k-j>\ell(\xi)+1$, then we must prove that
\begin{equation}\label{strict} ji_{k-1}+ (j+1)i_k+2(i_{k+1}+\cdots+i_\ell)< (\ell-k+j)\xi_{\ell-1}+\xi_\ell. \end{equation}
In other words, we must prove that if $\ell(\xi)+3\le k\le \ell+1$ and $1\le j<k-\ell(\xi)-1$, then
\begin{equation}\label{strict2} ji_{k-1}+ (j+1)i_k+2(i_{k+1}+\cdots+i_\ell)<  (\ell-k+j)\xi_{\ell-1}+\xi_\ell.\end{equation}
We prove this by a downward induction on $k$. If $k=\ell+1$ and $1\le j<\ell-\ell(\xi)$, then we have $$ji_\ell<j\xi_\ell\le (j-1)\xi_{\ell-1}+\xi_\ell,$$ and hence induction begins.

Assume now that we have proved the result for all $r>k$.  To prove the result for $k$, we proceed by an induction on $j$.  Taking  $j=1$ and $r=k+1$, we have, \begin{equation}\label{stepone}i_k+2(i_{k+1}+\cdots+i_\ell)< (\ell-k)\xi_{\ell-1}+\xi_\ell.\end{equation} Taking $j=2$ and $r=k$, (notice that this is allowed since $k-\ell(\xi)\ge 3$, i.e., $k\ge 3$), we have,
\begin{equation}\label{steptwo}2i_{k-1}+3i_k+2(i_{k+1}+\cdots+i_\ell)\le \sum_{p=k-2}^\ell\xi_p=(\ell-k+2)\xi_{\ell-1}+\xi_\ell.\end{equation}
Adding \eqref{stepone} and \eqref{steptwo} gives,
$$2(i_{k-1}+2i_k+2(i_{k+1}+\cdots+ i_\ell))<2(\ell-k+1)\xi_{\ell-1}+2\xi_\ell,$$\ \ i.e.,$$ i_{k-1}+2i_k+2(i_{k+1}+\cdots+ i_\ell)<(\ell-k+1)\xi_{\ell-1}+\xi_\ell.$$ This shows that \eqref{strict2} holds for $j=1$ and $k$. Assume that we have proved the result for all $1\le j'<j<k-\ell(\xi)-1$. Then taking $j'=j-1$ we have,
\begin{equation}\label{strict10} (j-1)i_{k-1}+ ji_k+2(i_{k+1}+\cdots+i_\ell)<(\ell-k+j-1)\xi_{\ell-1}+\xi_\ell.\end{equation}
On the other hand, since $2\le j<k-\ell(\xi)-1\le k-1$, we have $j+1\le k-1$ and so, we have the inequality
\begin{equation}\label{strict11} (j+1)i_{k-1}+ (j+2)i_k+2(i_{k+1}+\cdots+i_\ell)\leq (\ell-k+j+1)\xi_{\ell-1}+\xi_\ell. \end{equation}
Adding equations \ref{strict10} and \ref{strict11} gives,
$$2(ji_{k-1}+(j+1)i_k+2(i_{k+1}+\cdots+i_\ell))<2(\ell-k+j)\xi_{\ell-1}+2\xi_\ell,$$ which proves the inductive step for $j$ and completes the proof.
\end{pf}
\subsection{} The proof of Theorem \ref{sl2}(iii) is completed as follows. We proceed by induction on $\ell$. If $\ell=1$ then $V(\xi)\cong\ev_0 V(\xi_1)$. It is an elementary result in the representation theory of $\lie{sl}_2$ that $\{(x^-)^iv_{\xi_1}:0\le i\le \xi_1\}$ is a basis for $\ev_0 V(\xi_1)$. For the inductive step, note that if $(i_1,\cdots ,i_{\ell-1})\in\mathbb I(\xi^-)$, $$\varphi^-((x^-\otimes 1)^{i_1}\cdots (x^-\otimes t^{\ell-2})^{i_{\ell-1}}v_{\xi^-})=(x^-\otimes 1)^{i_1}\cdots (x^-\otimes t^{\ell-2})^{i_{\ell-1}}(x^-\otimes t^{\ell-1})^{\xi_\ell}v_{\xi}.$$ Since $\varphi^-$ is injective, it follows from the inductive hypothesis that the elements $$\left\{(x^-\otimes 1)^{i_1}\cdots (x^-\otimes t^{\ell-2})^{i_{\ell-1}}(x^-\otimes t^{\ell-1})^{\xi_\ell}v_{\xi}: (i_1,\cdots, i_{\ell-1},\xi_\ell)\in\mathbb I(\xi^-;\xi_\ell)\right\},$$ are a basis of the image of $\varphi^-$.

If $(i_1,\cdots, i_\ell)\in\mathbb I(\xi)$, with $i_\ell<\xi_\ell$, then by the inductive hypothesis, the image under $\varphi^+$ of the elements $$\left\{(x^-\otimes 1)^{i_1}\cdots (x^-\otimes t^{\ell-2})^{i_{\ell-1}}(x^-\otimes t^{\ell-1})^{i_\ell}v_{\xi}: (i_1,\cdots, i_{\ell-1},i_\ell)\in\mathbb I(\xi),\ \ i_\ell<\xi_\ell\right\},$$ is a basis for $V(\xi^+)$.

Theorem \ref{sl2} (iii) is now immediate from Theorem \ref{sl2} (i) and Proposition \ref{basissl2}.

\end{document}